\RequirePackage{ifpdf}
\ifpdf 
\documentclass[pdftex]{sigma}
\else
\documentclass{sigma}
\fi

\usepackage{algorithmic}

\newcommand{\cJ}{\cal{J}}

\newcommand{\cL}{\cal{L}}

\newcommand{\Q}{\mathbb{Q}}
\newcommand{\Z}{\mathbb{Z}}
\newcommand{\K}{\mathbb{K}}

\newcommand{\R}{\mathbb{R}}

\newcommand{\lcm}{\mathop{\mathrm{lcm}}\nolimits}
\newcommand{\lm}{\mathop{\mathrm{lm}}\nolimits}

\newcommand{\lc}{\mathop{\mathrm{lc}}\nolimits}
\newcommand{\pol}{\mathop{\mathrm{pol}}\nolimits}
\newcommand{\anc}{\mathop{\mathrm{anc}}\nolimits}

\newcommand{\Dp}{\mathop{\mathrm{dp}}\nolimits}
\newcommand{\Card}{\mathop{\mathrm{Card}}\nolimits}

\newcommand{\Id}{\mathop{\mathrm{Id}}\nolimits}
\newcommand{\idx}{\mathop{\mathrm{idx}}\nolimits}

\newenvironment{algorithm}[1]{
  \begin{center}
    {\bf Algorithm: #1} \\
     \begin{tabular}{|p{150mm}|} \hline
} {
 \\ \hline
 \end{tabular}
 \end{center}
}

\begin{document}
\allowdisplaybreaks

\renewcommand{\PaperNumber}{051}

\FirstPageHeading

\ShortArticleName{Generation of Dif\/ference Schemes for PDEs}

\ArticleName{Gr\"obner Bases and Generation of Dif\/ference
Schemes\\ for Partial Dif\/ferential Equations}

\Author{Vladimir P. GERDT~$^\dag$, Yuri A. BLINKOV~$^\ddag$ and
Vladimir V. MOZZHILKIN~$^\ddag$} \AuthorNameForHeading{V.P. Gerdt,
Yu.A. Blinkov and V.V. Mozzhilkin}

\Address{$^\dag$~Laboratory of Information Technologies,
           Joint Institute for Nuclear Research,\\
$\phantom{^\dag}$~141980 Dubna, Russia}

\EmailD{\href{mailto:gerdt@jinr.ru}{gerdt@jinr.ru}}

\URLaddressD{\url{http://compalg.jinr.ru/CAGroup/Gerdt/}}

\Address{$^\ddag$~Department of Mathematics and Mechanics, Saratov
University, 410071 Saratov, Russia}

\EmailD{\href{mailto:BlinkovUA@info.sgu.ru}{BlinkovUA@info.sgu.ru}}

\URLaddressD{\url{http://www.sgu.ru/faculties/mathematics/chairs/alg/blinkov.php}}

\ArticleDates{Received December 07, 2005, in f\/inal form April
24, 2006; Published online May 12, 2006}

\Abstract{In this paper we present an algorithmic approach to the
generation of fully conservative dif\/ference schemes for linear
partial dif\/ferential equations. The approach is based on
enlargement of the equations in their integral conservation law
form by extra integral relations between unknown functions and
their derivatives, and on discretization of the obtained system.
The structure of the discrete system depends on numerical
approximation methods for the integrals occurring in the enlarged
system. As a result of the discretization, a system of linear
polynomial dif\/ference equations is derived for the unknown
functions and their partial derivatives. A dif\/ference scheme is
constructed by elimination of all the partial derivatives. The
elimination can be achieved by selecting a proper elimination
ranking and by computing a Gr\"obner basis of the linear
dif\/ference ideal generated by the polynomials in the discrete
system. For these purposes we use the dif\/ference form of
Janet-like Gr\"obner bases and their implementation in Maple. As
illustration of the described methods and algorithms, we construct
a number of dif\/ference schemes for Burgers and Falkowich--Karman
equations and discuss their numerical properties.}


\Classification{68W30; 65M06; 13P10; 39A05; 65Q05}

\section{Introduction}

It is well-known that f\/inite dif\/ferences along with f\/inite
elements and f\/inite volumes are most important discretization
schemes for numerical solving of partial dif\/ferential equations
(PDEs)  (see, for
example,~\cite{GR'87,Str'89,GV'96a,QV'97,Th1,Th2,Sam'01,MM'05}).

Mathematical operations used in the construction of dif\/ference
schemes for PDEs are substantially symbolic. Thereby, it is a
challenge for computer algebra to provide an algorithmic tool for
automatization of the dif\/ference schemes constructing as well as
for the investigation of properties of the dif\/ference schemes.
One of the most fundamental requirements for a dif\/ference scheme
is its stability which can be analyzed with the use of computer
algebra methods and software~\cite{GV'96b}.

Furthermore, if PDEs admit a conservation law form or/and have
some symmetries, it is worthwhile to preserve these features at
the level of dif\/ference schemes too. In particular, a tool for
automatic construction of dif\/ference schemes should produce
conservative schemes whenever the original PDEs can be written in
the integral conservation law form. One of such tools GRIDOP
written in Reduce~\cite{LS'91,LSS'94} is based on symbolic
operator methods and generates conservative f\/inite-dif\/ference
schemes on rectangular domains in an arbitrary number of
independent variables. However, the generation is not entirely
automatic. A user of GRIDOP has to specify function spaces
together with associated scalar products and def\/ine grid
operators as f\/inite-dif\/ference schemes. Then the user may
provide partial dif\/ferential equations in terms of the def\/ined
grid operators or the adjoints of those operators. Under these
conditions the package returns the f\/inite-dif\/ference equations
for the dependent variables.

Besides, a few other applications of computer algebra are known to
construct f\/inite-dif\/ference schemes \cite{Fournie'99,Kol'01}
which, being also not completely automatic, are applicable to PDEs
of a certain form.

In this paper we describe a universal algorithmic approach to the
automatic generation of conservative dif\/ference schemes for
linear PDEs with two independent variables admitting the
conservation law form. This approach generalizes and extends the
observations of paper~\cite{MB'01} where it was noticed that a
conservative dif\/ference scheme can be derived as a compatibility
condition for a system of dif\/ference equations. The system is
composed of a discrete form of the original PDEs taken in the
integral conservation law form and of a number of natural integral
relations between functions and their partial derivatives. The
f\/inite-dif\/ference scheme is obtained by elimination of all the
partial derivatives from the system. We also show, by the example
of Burgers equation, that one can also apply the dif\/ference
elimination approach to generate of dif\/ference schemes without
use of conservation law form.

To perform the dif\/ference elimination we apply the Gr\"obner
bases method invented 40 years ago by Buchberger~\cite{Buch'65}
for polynomial ideals. This method has become the most universal
algorithmic tool in commutative algebra and algebraic geometry and
found also numerous fruitful applications for computations in
certain noncommutative polynomial rings as well as in rings of
linear dif\/ferential operators and dif\/ferential
polynomials~\cite{BW'98}. Nowadays, all modern general-purpose
computer algebra systems, for example, Maple~\cite{Maple} and
Mathematica~\cite{Math}, have special built-in modules
implementing algorithms for computing Gr\"obner bases. However,
the fastest implementation of these algorithms for commutative
polynomial algebra is done in the special-purpose systems
Singular~\cite{Singular} and Magma~\cite{Magma}. As to the
dif\/ference algebra~\cite{Cohn}, in spite of known for long time
(see~\cite{MLPK'99} and references therein) extensions of
Buchberger's algorithm~\cite{Buch'85} to dif\/ference polynomial
rings, there are only a few implementations of the algorithm
specialized to shift Ore algebra: in the Ore algebra library
package of Maple~\cite{Chyzak1}, in the library {\it
OreModules}~\cite{Chyzak2} developed using the latter package and
in Singular (Plural)~\cite{Plural}. These packages can be used
for computing Gr\"obner bases of linear dif\/ference ideals and
modules, and, in particular, for those linear systems which are
considered below.

In the given paper we present, however, another algorithm for
computing dif\/ference Gr\"obner bases. This algorithm is superior
over our Janet division algorithm whose polynomial
version~\cite{Gerdt'05} in most cases is computationally more
ef\/f\/icient than Buchberger's algorithm~\cite{ginv}. In
addition, unlike the above mentioned implementations of
Buchberger's algorithm for the shift Ore algebra, the algorithm
described below and its recent implementation in
Maple~\cite{GR'05} admit a natural extension to nonlinear
dif\/ference systems exactly in the same way as dif\/ferential
involutive algorithms~\cite{Janet,Thomas,Gerdt'99}. The algorithm
improves our Janet-like division algorithm~\cite{GB'05} adapted to
linear dif\/ference ideals~\cite{G'05}. The improvement includes,
in particular, the dif\/ference form of the involutive
criteria~\cite{Gerdt'05} modif\/ied for Janet-like reductions.
These criteria allow to avoid some useless reductions, and thereby
accelerate the computation.

The structure of the paper is as follows. In Section~2 we describe
the basic idea of our approach to the generation of
f\/inite-dif\/ference schemes for two-dimensional PDEs. Section~3
contains def\/initions and notions of dif\/ference algebra which
are used in the sequel. Here we def\/ine Gr\"obner bases for
linear dif\/ference ideals and their special form called
Janet-like bases. In Section~4 we present an improved version of
the algorithm in paper~\cite{G'05} and brief\/ly discuss some
relevant computational aspects. Section~5 illustrates our approach
to construction of dif\/ference schemes by simplest second-order
equations -- Laplace's equation, the wave equation, the heat
equation, and by the f\/irst-order advection equation. In
Section~6 we generate several dif\/ference schemes for Burgers
equation. But for all that, to avoid problems arising in computing
of nonlinear Gr\"obner bases, we denote the square of the
dependent variable by an extra function. Besides, we characterize
some of the constructed schemes by the modif\/ied equation method.
In Section~7 we consider the two-dimensional quadratically
nonlinear Falkowich--Karman equation describing transonic f\/low
in gas dynamics. Here, we succeeded in computing of the nonlinear
Gr\"obner basis by hand, and in that way generated the cubic
nonlinear dif\/ference scheme which possesses some attractive
properties. These properties as well as those of the schemes
generated for Burgers equation are illustrated by some numerical
experiments in Section~8. We conclude in Section~9.

\section{Basic idea}

It is well-known~\cite{QV'97,Th1,Th2,MM'05} that a rather wide
class of scalar PDE and some systems of PDEs can be written in the
conservation law form
\begin{gather}
\frac{\partial \boldsymbol{v}}{\partial x} +
\frac{\partial}{\partial y}\boldsymbol{F}(\boldsymbol{v})=0,
\label{cons_law}
\end{gather}
where $\boldsymbol{v}$ is a $m$-vector function in the unknown
$n$-vector function $\boldsymbol{u}$ and its partial derivatives
$\boldsymbol{u}_x,\boldsymbol{u}_y$, $\boldsymbol{u}_{xx},\ldots$.
The vector function $\boldsymbol{F}$ maps $R^m$ into $R^m$.

By Green's theorem (curl theorem in the plane), vector
PDE~(\ref{cons_law}) is equivalent to the integral relation
\begin{gather}
 \oint _{\Gamma} - \boldsymbol{F}(\boldsymbol{v}) dx + \boldsymbol{v}dy = 0, \label{int_cons_law}
\end{gather}
where $\Gamma$ is an arbitrary closed contour. Approximation
of~(\ref{int_cons_law}) rather than of~(\ref{cons_law}) on a
dif\/ference grid (balance or integro-interpolation method) is a
natural way to generate conservative f\/inite-dif\/ference schemes
for PDEs of order two and higher.

Throughout this paper we shall consider orthogonal and uniform
grids with the grid mesh steps $h_1$ and $h_2$
\begin{gather}
x_{j+1} -  x_{j} = h_1,\qquad y_{k+1} -  y_{k} = h_2
\label{discret}
\end{gather}
and denote the grid values of the vector function
$\boldsymbol{u}(x, y)$ and all its partial derivatives occurring
in~(\ref{int_cons_law}) by
\begin{gather}
\boldsymbol{u}(x, y) \Longrightarrow \boldsymbol{u}(x_j,y_k)\equiv
\boldsymbol{u}_{jk},\qquad
{\boldsymbol{u}}_x(x, y) \Longrightarrow {\boldsymbol{u}}_x(x_j,y_k)\equiv ({\boldsymbol{u}}_x)_{jk}, \nonumber\\
{\boldsymbol{u}}_y(x, y) \Longrightarrow
{\boldsymbol{u}}_y(x_j,y_k)\equiv ({\boldsymbol{u}}_y)_{jk},\qquad
 {\boldsymbol{u}}_{xx}(x, y) \Longrightarrow {\boldsymbol{u}}_{xx}(x_j,y_k)\equiv ({\boldsymbol{u}}_{xx})_{jk}, \label{grid} \\
\cdots\cdots\cdots\cdots\cdots\cdots\cdots\cdots\cdots\cdots\cdots\cdots\cdots\cdots\cdots\cdots\cdots\cdots
\cdots\cdots\cdots\cdots\cdots\cdots\nonumber
\end{gather}

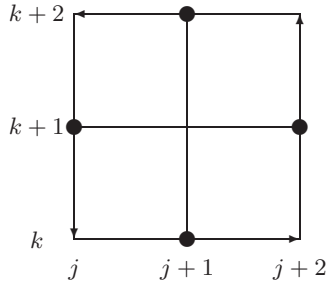
\begin{figure}[t]
\unitlength=1.00mm \special{em:linewidth 0.4pt}
\linethickness{0.4pt}
\begin{center}
\begin{picture}(110.00,36.00)(-18, 0)
\put(21.00,5.00){\vector(1,0){30.00}}
\put(51.00,5.00){\vector(0,1){30.00}}
\put(51.00,35.00){\vector(-1,0){30.00}}
\put(21.00,35.00){\vector(0,-1){30.00}}
\put(21.00,20.00){\line(1,0){30.00}}
\put(36.00,5.00){\line(0,1){30.00}}
\put(36.00,35.00){\circle*{2.00}}
\put(51.00,20.00){\circle*{2.00}}
\put(21.00,20.00){\circle*{2.00}} \put(36.00,5.00){\circle*{2.00}}
\put(16.00,5.00){\makebox(0,0)[cc]{\footnotesize $k$}}
\put(16.00,20.00){\makebox(0,0)[cc]{\footnotesize $k+1$}}
\put(16.00,35.00){\makebox(0,0)[cc]{\footnotesize $k+2$}}
\put(21.00,1.00){\makebox(0,0)[cc]{\footnotesize $j$}}
\put(36.00,1.00){\makebox(0,0)[cc]{\footnotesize $j+1$}}
\put(51.00,1.00){\makebox(0,0)[cc]{\footnotesize $j+2$}}
\end{picture}
\vspace{-2mm} \caption{Integration contour on
grid.}\label{int_cont}
\end{center}\vspace{-7mm}
\end{figure}

Choose the integration contour as follows (see
Fig.~\ref{int_cont}) and add all the related integral relations
between the dependent vector variable and its partial derivatives:
\begin{gather}
\int_{x_j}^{x_{j+2}}  \boldsymbol{u}_x dx =
\boldsymbol{u}(x_{j+2}, y) - \boldsymbol{u}(x_{j}, y),\qquad
\int_{y_k}^{y_{k+2}}  \boldsymbol{u}_y dy = \boldsymbol{u}(x,
y_{k+2}) -
\boldsymbol{u}(x, y_{k}), \nonumber\\
 \int_{x_j}^{x_{j+2}}  \boldsymbol{u}_{xx} dx = \boldsymbol{u}_x(x_{j+2}, y) -
\boldsymbol{u}_x(x_{j}, y),\qquad \int _{y_k}^{y_{k+2}}
\boldsymbol{u}_{xy} dy = \boldsymbol{u}_x(x, y_{k+2}) -
\boldsymbol{u}_x(x, y_{k}),  \label{rel}\\
\cdots\cdots\cdots\cdots\cdots\cdots\cdots\cdots\cdots\cdots\cdots\cdots\cdots\cdots\cdots\cdots\cdots\cdots
\cdots\cdots\cdots\cdots\cdots\cdots\cdots\cdots\cdots\nonumber
\end{gather}

\newpage

To obtain a source system of discrete equations for constructing
of a dif\/ference scheme, we consider numerical approximations of
the integral equations~(\ref{int_cons_law}) for the contour of
Fig.~\ref{int_cont} and of the relations~(\ref{rel}) in terms of
the grid unknowns~(\ref{grid}). Although generally one can use
dif\/ferent numerical approximations for the integral equations
in~(\ref{int_cons_law}) and (\ref{rel}), we apply here for all
these equations the simplest rectangle (midpoint) rule
\begin{gather}
  \boldsymbol{F}(\boldsymbol{v})_{j+1 \, k+2}  -\boldsymbol{F}(\boldsymbol{v})_{j+1 \, k} +
  (\boldsymbol{v}_{j+2 \, k+1} - \boldsymbol{v}_{j \, k+1})   = 0, \nonumber\\
 (\boldsymbol{u}_x)_{j+1 \, k} \cdot 2 h_1 = \boldsymbol{u}_{j+2 \, k} - \boldsymbol{u}_{j \, k},\qquad
 (\boldsymbol{u}_y)_{j \, k+1} \cdot 2 h_2 = \boldsymbol{u}_{j \, k+2} - \boldsymbol{u}_{j \, k},\label{disc_system}\\
\cdots\cdots\cdots\cdots\cdots\cdots\cdots\cdots\cdots\cdots\cdots\cdots\cdots\cdots\cdots\cdots\cdots\cdots\cdots\cdots
\cdots\cdots \nonumber
\end{gather}

Thereby, we obtained system~(\ref{disc_system}) of dif\/ference
equations in grid unknowns~(\ref{grid}). In doing so, the number
of scalar equations added to the (vector) integral
equation~(\ref{int_cons_law}) corresponds to the number of proper
partial derivatives of order less than or equal to the order of
partial derivatives involved in the integrand
of~(\ref{int_cons_law}).

It follows that eliminating from~(\ref{disc_system}) all the
proper grid partial derivatives gives equations containing only
independent (vector) function $\boldsymbol{u}$, and, hence
composing a f\/inite-dif\/ference scheme. If
system~(\ref{disc_system}) is linear, then this dif\/ference
elimination can always be algorithmically achieved by the
Gr\"obner bases method considered in the next section.

It should be noted that generation of f\/inite-dif\/ference
schemes on grid~(\ref{discret}) by the elimination can be also
applied to PDEs irrelative to their conservation law properties.
Again, one has to add to the initial dif\/ferential equations,
written in terms of grid variables (\ref{grid}), the corresponding
number of integral relations~(\ref{rel}) and approximate them by
numerical quadrature formulas. Such an approach may give more
f\/lexibility in generation of distinct dif\/ference schemes, and
we apply it in Section~4 to the f\/irst-order advection equation,
and in Section~5 to Burgers equation. Gene\-ral\-ly, however, for
the second-(and higher-) order PDEs admitting the integral
conservation law form, the dif\/ference scheme obtained directly
from the dif\/ferential form may not be conservative. Besides, the
dif\/ference elimination based on the integral form is usually
more ef\/f\/icient than that based on the dif\/ferential form.
This is because the number of partial derivatives to be eliminated
in the former case is smaller than in the latter case. Indeed, the
integrand in~(\ref{int_cons_law}) has the dif\/ferential order
smaller by one than that in~(\ref{cons_law}) whereas computational
complexity of the elimination is at least exponential in the
number of eliminated variables~\cite{ZuGathen}.

\section[Difference Gr\"obner bases]{Dif\/ference Gr\"obner bases}

A {\it difference ring} $R$ is a commutative ring with a unity
together with a f\/inite set of mutually commuting injective
endomorphisms $\theta_1,\ldots,\theta_n$ of $R$. Similarly, one
def\/ines a {\it difference field}. Elements $\{y^1,\ldots,y^m\}$
in a dif\/ference ring containing $R$ are said to be {\it
difference indeterminates} over $R$ if the set
\[
\big\{\theta_1^{k_1}\cdots \theta_1^{k_1}y^j\mid
\{k_1,\ldots,k_n\}\in \Z^n_{\geq 0},\ 1\leq j\leq m\big\}
\]
is algebraically independent over $R$.
\newpage

Hereafter we shall consider the ring of functions of $n$ variables
$x_1,\ldots,x_n$ with the basis endomorphisms $ \theta_i\circ
f(x_1,\ldots,x_n)=f(x_1,\ldots,x_i+1,\ldots,x_n). $ acting as
shift operators.

The f\/ield $\Q(x_1,\ldots,x_n)$ of rational functions in
$\{x_1,\ldots,x_n\}$ whose coef\/f\/icients are rational numbers
is an example of dif\/ference f\/ield, and we shall assume in the
next sections that the coef\/f\/icients of PDEs belong to this
f\/ield.

Let $\K$ be a dif\/ference f\/ield, and $\R:=\K\{y^1,\ldots,y^m\}$
be the dif\/ference ring of polynomials over $\K$ in variables
$\{\theta^\mu \circ y^k \mid \mu\in \Z^n_{\geq 0},\,k=1,\ldots,m
\}$. Hereafter, we denote by $\R_L$ the set of linear polynomials
in $\R$ and use the notations:
\begin{gather}
\Theta:=\{\theta^\mu \mid \mu\in \Z^n_{\geq 0} \},\qquad
\deg_i(\theta^\mu\circ y^k):=\mu_i,  \qquad
\deg(\theta^\mu \circ y^k):=|\mu|:=\sum_{i=1}^n \mu_i, \nonumber\\
\lcm(\mu,\nu):=\{\max\{\mu_1,\nu_1\},\ldots,\max\{\mu_n,\nu_n\}\},\qquad
\lcm(\theta^{\mu}\circ y^k,\theta^{\nu}\circ y^k):=\theta^{\lcm(\mu,\nu)}\circ y^k,\nonumber\\
\theta^\mu\circ y^k\, \sqsubset\, \theta^\nu \circ y^k\ \
\mathrm{when}\ \ \nu-\mu\in \Z^n_{\geq 0}\ \wedge\
|\nu-\mu|>0.\label{theta}
\end{gather}
A {\it difference ideal} is an ideal $I \subseteq \R$ closed under
the action of any operator from $\Theta$. If
$F:=\{f_1,\ldots,f_k\}\subset \R$ is a f\/inite set, then the
smallest dif\/ference ideal containing $F$ will be denoted by
$\Id(F)$. If for an ideal $I$ there is $F\subset \R_L$ such that
$I=\Id(F)$, then $I$ is a {\it linear difference ideal}.

A total ordering $\succ$ on the set of $\theta^\mu \circ y^{\,j}$
is a {\it ranking} if $\forall \,i,j,k,\mu,\nu$ the following
hold:
\[
\theta_{i} {\theta^\mu \circ y^{\,j}} \succ {\theta^\mu}\circ
y^{\,j},\qquad \theta^\mu \circ y^{\,j} \succ \theta^\nu \circ y^k
\iff {\theta_i}  {\theta^\mu \circ y^{\,j}}
   \succ {\theta_i}{\theta^\nu} \circ y^k.
\]

If $|\mu| \succ |\nu| \Longrightarrow {\theta^\mu \circ y^{\,j}}
\succ {\theta^\nu} \circ y^k$ the ranking is {\it orderly}. If $j
> k \Longrightarrow {\theta^\mu} \circ y^{\,j} \succ {\theta^\nu}
\circ y^k$ the ranking is {\it elimination}.

Given a ranking $\succ$, a linear polynomial $f\in \R_L\setminus
\{0\}$ has the {\it leading term} of the form $a\,\vartheta \circ
y^j$, $\vartheta\in \Theta$, where $\vartheta \circ y^j$ is
maximal w.r.t. $\succ$ among all $\theta^\mu \circ y^k$ which
appear with nonzero coef\/f\/icient in $f$. $\lc(f):=a\in
\K\setminus \{0\}$ is the {\it leading coefficient} and
$\lm(f):=\vartheta \circ y^{\,j}$ is the {\it leading (head)
monomial}.

A ranking acts in $\R_L$ as a {\it monomial order}. If $F
\subseteq \R_L \setminus \{ 0 \}$, then $\lm(F)$ will denote the
set of the leading monomials and $\lm_j(F)$ will denote its subset
for indeterminate $y^{\,j}$. Thus,
\[
\lm(F)=\cup_{j=1}^m \lm_j(F).
\]

Given a nonzero linear dif\/ference ideal $I=\Id(G)$ and a ranking
$\succ$, the ideal generating set $G=\{g_1,\ldots,g_s\}\subset
\R_L$ is a {\it Gr\"obner basis}~\cite{MLPK'99,Buch'85} of $I$ if
\begin{gather}
\forall \, f\in I\cap \R_L\setminus \{0\},\quad \exists \, g\in
G,\ \theta \in \Theta\  :\quad \lm(f)=\theta \circ \lm(g).
\label{GB}
\end{gather}
It follows that the head monomial of $f\in I\setminus \{0\}$, as
well as the polynomial $f$ itself, is {\it reducible modulo $G$}
and yields {\it the head reduction}:
\[
f \xrightarrow[g]{} f':=f-\lc(f)\,\theta \circ (g/\lc(g)),\qquad
f'\in I.
\]
If $f'\neq 0$, then its leading monomial is again reducible modulo
$G$, and, by repeating the reduction f\/initely many
times~\cite{BW'98,MLPK'99,Buch'85} we obtain $ f
\xrightarrow[G]{}0$. Generally, if a polynomial $h\in \R_L$
contains a term with monomial $u$ and coef\/f\/icient $c\neq 0$
such that $u=\vartheta \circ \lm(f)$ for some $\vartheta \in
\Theta$ and $f\in F\subset \R_L\setminus \{0\}$, then $h$ can be
reduced:
\begin{gather}
h \xrightarrow[g]{} h':=h-c\,\theta \circ (f/\lc(f)).
\label{elem_red}
\end{gather}
By applying the reduction f\/initely many times, one obtains a
polynomial $\bar{h}$ which is either zero or such that all its
(nonzero) terms have monomials {\it irreducible modulo the set $F
\subset \R_L$}. In both cases $\bar{h}$ is said to be in the {\it
normal form modulo $F$} (denotation: $\bar{h}=NF(h,F)$). A
Gr\"obner basis $G$ is {\it reduced} if
\begin{displaymath}
\forall \, g\in G\ :\quad  g=NF(g,G\setminus \{g\}).
\end{displaymath}

In our algorithmic construction of reduced Gr\"obner bases we
shall use a restricted set of reductions called {\it Janet-like}
(cf.~\cite{GB'05,G'05}) and def\/ined as follows.

For a f\/inite set $F \subseteq \R_L$ and a ranking $\succ$, we
partition every $\lm_k(F)$ into groups labeled by
$d_0,\ldots,d_i\in \Z_{\geq 0}$,\ $(0\leq i\leq n)$. Here
$[0]_k:=\lm_k(F)$ and for $i>0$ the group $[d_0,\ldots,d_i]_k$ is
def\/ined as
\[
[d_0,\ldots,d_i]_k:=\{u\in \lm_k(F) \mid d_0=0,d_j=\deg_j(u),\
1\leq j\leq i \}.
\]
Now we characterize a monomial $u\in \lm_k(F)$ by the nonnegative
integer $\lambda_i$:
\[
\lambda_i(u,\lm_k(F)):=\max\{\deg_i(v) \mid u,v\in
[d_0,\ldots,d_{i-1}]_k\}-\deg_i(u).
\]
If $\lambda_i(u,\lm_k(F))>0$, then $\theta_i^{s_i}$ such that
\[
s_i:=\min\{\deg_i(v)-\deg_i(u) \mid u,v\in
[d_0,\ldots,d_{i-1}]_k,\, \deg_i(v)>\deg_i(u)\}
\]
is called a {\it difference power} for $f\in F$ with $\lm(f)=u$.

Let $DP(f,F)$ denotes the set of all dif\/ference powers for $f\in
F$. Now we def\/ine the partition of the set $\Theta$ into two
disjoint subsets
\[
\bar{\Theta}(f,F):=\{\theta^\mu \mid \exists \, \theta^\nu \in
DP(f,F)\ :\ \mu-\nu\in \Z^n_{\geq 0} \},\qquad {\cJ}(f,F):=\Theta
\setminus \bar{\Theta}(f,F),
\]
which is similar to the partition of monomials into
nonmultiplicative and multiplicative ones in the involutive
approach~\cite{Gerdt'05}.

A f\/inite basis $G$ of $I=\Id(G)$ is called
Janet-like~\cite{GB'05,G'05} if
\begin{gather}
\forall \, f\in I\cap \R_L\setminus \{0\}, \ \exists \, g\in G,
\theta \in {\cJ}(g,G)\  :\quad  \lm(f)=\theta \circ
  \lm(g). \label{JLGB}
\end{gather}
In full analogy with~(\ref{elem_red}) a {\it ${\cJ}$-reduction} is
def\/ined as
\begin{gather}
h \xrightarrow[g]{} h':=h-c\,\theta \circ (f/\lc(f)),\qquad
\theta\in {\cJ}(f,F), \label{elem_J_red}
\end{gather}
for polynomial $h\in \R_L$ containing monomial $u$ with
coef\/f\/icient $c\neq 0$ satisfying $u=\vartheta \circ \lm(f)$
for some $f\in F\subset \R_L\setminus \{0\}$  and $\vartheta \in
{\cJ}(f,F)$.

Apparently, any element in the ideal $I=\Id(G)$ is {\it $\cJ$-head
reduced} to zero by the f\/inite sequence of {$\cJ$-head
reductions} by elements $g\in G$ in the Janet-like basis $G$:
\begin{gather}
f \xrightarrow[g]{} f':=f-\lc(f)\,\theta \circ (g/\lc(g)),\qquad
\theta\in {\cJ}(g,G),\qquad f'\in I. \label{elem_J_head_red}
\end{gather}

If the leading monomial of $p\in \R\setminus \{0\}$ is not
$\cJ$-reducible modulo a f\/inite subset $F\subset \R\setminus
\{0\}$ we say that $p$ is in the {$\cJ$-head normal form modulo
$F$} and write $p=HNF_{\cJ}(p,F)$. If none of monomials in $p$ is
$\cJ$-reducible modulo $F$ we say that $p$ is in {\it the (full)
normal form modulo $F$} and write $p=NF_{\cJ}(p,F)$.

Since $\cJ$-reducibility implies the Gr\"obner reducibility
(\ref{elem_red}), a Janet-like basis satisfying~(\ref{JLGB}) is
a~Gr\"obner basis. The converse is generally not true, that is,
not every Gr\"obner basis is Janet-like.

The properties of a Janet-like basis are very similar to those of
a Janet basis~\cite{GB'98}, but the former is generally more
compact than the latter. For all that we consider hereafter only
minimal bases.

Let $GB$ be a reduced Gr\"obner basis, satisfying~\cite{Buch'85}:
\begin{gather}
\forall\, g\in GB\ :\quad g=NF(g,GB\setminus \{g\}).
\label{red_GB}
\end{gather}
Let now $JB$ be a Janet basis, and $JLB$ be a Janet-like basis of
the same ideal and for the same ranking. Then for their
cardinalities the inclusion $\Card(GB)\leq \Card(JLB)\leq
\Card(JB)$ holds~\cite{Gerdt'05,GB'05}. Here $\Card$ abbreviates
{\it cardinality}, that is, the number of elements.

Whereas the algorithmic characterization of a Gr\"obner basis is
zero redundancy of all its $S$-polynomials~\cite{Buch'65,Buch'85},
the algorithmic characterization of a Janet-like basis $G$ has the
following form~(cf.~\cite{GB'05}):
\begin{gather}
\forall \, g\in G,\ \vartheta\in DP(g,G):\quad NF_{\cJ}(\vartheta
\circ g,G)=0. \label{alg_char}
\end{gather}
These conditions are at the root of the algorithmic construction
of Janet-like bases as described in the following section.

\section{Algorithm}

In this section we present an algorithm for constructing a reduced
Gr\"obner basis~(\ref{GB}) of the ideal generated by an input set
of linear dif\/ference polynomials. The algorithm is an improved
version of the algorithm in paper~\cite{G'05} and translates to
the dif\/ference form of the polynomial involutive
algorithm~\cite{Gerdt'05} modif\/ied for the Janet-like
reductions.

To apply the dif\/ference form of criteria to avoid some
unnecessary reductions we need the following def\/inition.

An {\it ancestor} of a dif\/ference polynomial $f\in F\subset
\R_L\setminus \{0\}$ is a polynomial $g\in F$ of the smallest
$\deg(\lm(g))$ among those satisfying $f=\theta \circ g$ modulo
$\Id(F\setminus \{f\})$ with $\theta\in \Theta$. If for all that
\begin{displaymath}
\deg(\lm(g))<\deg(\lm(f)),
\end{displaymath}
then the ancestor $g$ of $f$ is called {\it proper}.

If an intermediate polynomial $h$ that arose in the course of the
below algorithm has a proper ancestor $g$ in the intermediate
basis $G$, then $h$ has been obtained from $g$ via a sequence of
shift operations of the form $\vartheta \circ g$ where $\vartheta
\in DP(g,G)$ with $\lm(\vartheta \circ g)$ $\cJ$-irreducible
modulo $G$. For the ancestor $g$ itself the equality
$\lm(\anc(g))=\lm(g)$ holds.

In the main algorithm {\bf Gr\"{o}bnerBasis} and its subalgorithms
we endow every element $f\in G$ in the intermediate set of
dif\/ference polynomials $G$ (occuring in the set $T$) with a
triple structure of the form:
\begin{displaymath}
p=\{f,\, g,\, dpow\}\,,
\end{displaymath}
where
\begin{gather*}
\pol(p):=f\ \ \mbox{is the polynomial}\ f\ \mbox{itself},\\
\anc(p):=g\ \ \mbox{is an ancestor of}\ f\
 \mbox{in}\ G,\\
\Dp(p):=dpow\ \ \mbox{is a (possibly empty) subset of $DP(f,G)$}.
\end{gather*}

The set $dpow$ associated with the polynomial $f$ accumulates all
the dif\/ference powers for $f$ which have been already applied to
$f$ in the course of the algorithm. Keeping this information
serves to avoid useless repeated applications of the dif\/ference
power operators. Knowledge of ancestors for elements in the
intermediate basis  helps to avoid some unnecessary reductions by
applying Buchberger's chain criterion~\cite{Buch'85} adapted to
Janet-like reductions.

\begin{algorithm}{Gr\"{o}bnerBasis($F,\prec$)\label{Algorithm GB}}
\begin{algorithmic}[1]
\INPUT $F\in \R_L\setminus \{0\}$, a f\/inite set; $\prec$, a ranking\\
\OUTPUT $G$, a reduced Gr\"obner basis of $\Id(F)$ \STATE {\bf
choose} $f\in F$ with the lowest $\lm(f)$ w.r.t. $\succ$ \STATE
$T:=\{f,f,\varnothing\}$ \STATE $Q:=\{\, \{q,q,\varnothing\} \mid
q\in F\setminus \{f\}\, \}$ \STATE $Q:=${\bf
HeadReduce}$(Q,T,\succ)$
  \WHILE{$Q \neq \varnothing$}
    \STATE {\bf choose} $p \in Q$ with the lowest $\lm(\pol(p))$ w.r.t. $\succ$
    \STATE $Q:=Q\setminus \{p\}$
    \IF{$\pol(p) = \anc(p)$}
          \FORALL{$\{\ q \in T \mid \lm(\pol(q))=\theta^{\mu} \circ \lm(\pol(p)),\ |\mu|>0\ \}$}
          \STATE $Q:=Q \cup \{q\}$; \hspace*{0.3cm} $T:=T \setminus \{q\}$
          \ENDFOR
    \ENDIF
    \STATE $h:=${\bf{TailNormalForm}}($p,T,\prec)$
    \STATE $T:=T \cup \{h,\anc(p),\Dp(p)\}$
    \FORALL{$q\in T$ {\bf and} $\vartheta \in DP(q,T)\setminus \Dp(q)$}
      \STATE $Q:=Q \cup \{\{\vartheta \circ \pol(q),\anc(q),\varnothing\}\}$
      \STATE $\Dp(q):=\Dp(q)\cap DP(q,T)\cup \{\vartheta\}$
    \ENDFOR
    \STATE $Q:=${\bf HeadReduce}$(Q,T,\prec)$
  \ENDWHILE
  \RETURN $\{\pol(f)\mid f\in T\}$ or $\{\pol(f)\mid f\in T\mid f=\anc(f)\}$
\end{algorithmic}
\end{algorithm}

In the above main algorithm {\bf Gr\"{o}bnerBasis} and its
subalgorithms presented below, where no confusion can arise, we
simply refer to the triple set $T$ as the second argument in $DP$,
$NF_{\cJ}$, and $HNF_{\cJ}$ instead of the polynomial set
$\{g=\pol(t) \mid t\in T\}$. Sometimes we also refer to the triple
$p$ instead of $\pol(p)$. Besides, when we speak of reduction of
the triple set $Q$ modulo triple set $T$ we mean reduction of the
polynomial set
\[
\{f=\pol(q) \mid q\in Q\}
\]
modulo
\[
\{g=\pol(t) \mid t\in T\}.
\]

Correctness and termination of algorithm {\bf Gr\"{o}bnerBasis}
can be shown exactly as in the polynomial
case~\cite{Gerdt'05,GB'05}. Here we only elucidate some related
features of the algorithm.

At steps 4 and 19 the ${\cJ}$-head reduction is performed for the
dif\/ference polynomials in $Q$ modulo those in $T$. Then the
remaining tail reduction is done in line 13 to obtain the (full)
${\cJ}$-normal form. Thereby, the main {\bf while-}loop 5--20
terminates when the conditions~(\ref{alg_char}) hold for the
dif\/ference polynomial set $G$ composed from the f\/irst elements
of triples in $T$
\begin{gather}
G:=\{\pol(g)\mid g\in T\}, \label{G_in_T}
\end{gather}
and the set $Q$ is empty. The upper {\bf for}-loop 9--11 provides
minimality of the output Janet-like basis contained in
$T$~\cite{Gerdt'05}. Another {\bf for}-loop 15--18 constructs new
conditions~(\ref{alg_char}) which have to be further examined
because of the insertion of a new element in $T$ at step 14.
Besides, the set of dif\/ference powers is upgraded at step 17.

Furthermore, the main algorithm {\bf Gr\"{o}bnerBasis} together
with its subalgorithms presented below ensures that every element
in the output Janet-like basis composed from the f\/irst elements
in the triple set $T$ has one and only one ancestor. This ancestor
is apparently irreducible, in the Gr\"obner
sense~(\ref{elem_red}), by other elements in the basis. Thereby,
those elements in the output basis that have no proper ancestors
constitute the reduced Gr\"obner basis~(\ref{red_GB}) that is
returned by the main algorithm at the last step 21.

The algorithm {\bf HeadReduce} invoked in lines 4 and 19 of the
main algorithm returns the set~$Q$ which, if nonempty, contains
part of the intermediate basis ${\cJ}$-head reduced modulo $T$.
The reductions are performed by its subalgorithm {\bf
HeadNormalForm} that is invoked at step~6 of the algorithm.

If algorithm {\bf HeadNormalForm} returns $h\neq 0$, then $\lm(h)$
does not belong to the initial ideal generated by $\{\lm(\pol(f))
\mid f\in Q\cup T\}$~\cite{Gerdt'05,GB'98}. In this case the
triple $\{h,h,\varnothing\}$ for $h$ is inserted (line 9) into the
output set $Q$. Otherwise, the output set $Q$ retains the triple
$p$ as it is in the input.

In the case when $h=0$ and $\pol(p)$ has no proper ancestors that
is checked at step 14, all the descendant triples for $p$, if any,
are deleted from the intermediate set $S$ at step 16. Such
descendants cannot occur in $T$ owing to the choice conditions at
steps~1,~6  and to the displacement condition of step 9 in the
main algorithm {\bf Gr\"{o}bnerBasis}. Steps 14--18 serve for the
memory saving and can be ignored if the memory restrictions are
not very critical for a given problem. In this case all those
descendants will be casted away by the criteria checked in the
below algorithm {\bf HeadNormalForm}.

Algorithm {\bf HeadNormalForm} performs verif\/ication (step 3) of
${\cJ}$-head reducibility of the input polynomial $h$ modulo the
polynomial set~(\ref{G_in_T}). This verif\/ication consists in
searching a~dif\/ference polynomial (reductor) in the set $G$
def\/ined in~(\ref{G_in_T}) such that $G$ yields the
reduction~(\ref{elem_J_red}). If the search fails, that is, there
is no $\cJ$-reductor, the algorithm returns at step 4 the input
polynomial.

For the ${\cJ}$-head reducible input polynomial $\pol(p)$ that is
checked at step 3 of algorithm {\bf HeadNormalForm}, the following
three criteria are verif\/ied at step 9
\begin{gather}
Criteria(p,g)=C_1(p,g) \vee  C_2(p,g) \vee C_3(p,g)\,,
\label{crit}
\end{gather}
\noindent where
\begin{itemize}
\itemsep=0pt \item[] $C_1(p,g)$ is true $\Longleftrightarrow$
$\lcm(\lm(\anc(p)), \lm(\anc(g))) \sqsubset \lm(\pol(p)) $,

\item[] $C_2(p,g)$ is true $\Longleftrightarrow$ $\exists$ $t\in T$ such that \\
\hspace*{1.0cm}$\lcm(\lm(\pol(t)),\lm(\anc(p))) \sqsubset
\lcm(\lm(\anc(p)), \lm(\anc(g)))$
$\wedge$ \\
\hspace*{1.0cm}$\lcm(\lm(\pol(t)),\lm(\anc(g))) \sqsubset
\lcm(\lm(\anc(p)), \lm(\anc(g)))$,

\item[] $C_3(p,g)$ is true $\Longleftrightarrow$ $\exists$ $t\in
T$ $\wedge$
$y\in NM_{\cL}(t,T)$ with $\lm(\pol(t))\cdot y=\lm(\pol(p))$,  \\
\hspace*{1.0cm}$\lcm(\lm(\anc(p)),\lm(\anc(t))) \sqsubset
\lm(\pol(p))$
$\wedge$ $\idx(t,T)<\idx(f,T)$, \\
\hspace*{1.0cm}where $\idx(t,T)$ enumerates the position of triple
$t$ in set $T$.

\end{itemize}

\noindent In aggregate, criteria~(\ref{crit}) translate
(cf.~\cite{Gerdt'05,AH'02}) Buchberger's chain
criterion~\cite{Buch'85} into the linear dif\/ference algebra.

In addition, if all dif\/ference polynomials in the input set $F$
for the main algorithm {\bf Gr\"{o}bner\-Basis} have constant
coef\/f\/icients, then the set of criteria~(\ref{crit}) can be
enlarged with one more criterion~$C_4$:

\begin{itemize}
\itemsep=0pt \item[] $C_4(p,g)$ is true for $\lm(\pol(p))=\theta
\circ y^k$, $\lm(\pol(g))=\vartheta \circ y^k$
$\Longleftrightarrow$ $\lcm(\theta,\vartheta)=\theta \,
\vartheta$\,.
\end{itemize}
Criterion $C_4$ is the dif\/ference form (cf.~\cite{Gerdt'05}) of
Buchberger's co-prime criterion~\cite{Buch'85}.

\newpage
\noindent
\begin{algorithm}{HeadReduce$(Q,T,\prec)$}
\begin{algorithmic}[1]
\INPUT $Q$ and $T$, sets of triples; $\prec$, a ranking
\OUTPUT ${\cJ}$-head reduced set $Q$ modulo $T$
 \STATE $S:=Q$
 \STATE $Q:=\varnothing$
 \WHILE {$S \neq \varnothing$}
   \STATE {\bf choose} $p\in S$
   \STATE $S:=S\setminus \{p\}$
   \STATE $h:=${\bf{HeadNormalForm}}$(p,T)$
   \IF{$h\neq 0$}
       \IF{$\lm(\pol(p))\neq \lm(h)$}
         \STATE $Q:=Q\cup \{h,h,\varnothing\}$
       \ELSE
         \STATE $Q:=Q\cup \{p\}$
       \ENDIF
   \ELSE
      \IF{$\lm(\pol(p))=\lm(\anc(p))$}
         \FORALL{$\{ q\in S \mid \anc(q)=\pol(p)\}$}
          \STATE $S:=S\setminus \{q\}$
         \ENDFOR
      \ENDIF
   \ENDIF
\ENDWHILE \RETURN $Q$
\end{algorithmic}
\end{algorithm}

\noindent
\begin{algorithm}{HeadNormalForm($p,T,\prec$)}
\begin{algorithmic}[1]
\INPUT $T$, a set of triples; $p$, a triple;
 $\prec$, a ranking
\OUTPUT $h=HNF_{\cJ}(p,T)$, the ${\cJ}$-head normal form
   of $\pol(p)$ modulo $T$
   \STATE $h:=\pol(p)$
   \STATE $G:=\{\,\pol(g)\mid g\in T\,\}$
    \IF{$\lm(h)$ is ${\cJ}$-irreducible modulo $G$}
       \RETURN $h$
    \ELSE
       \STATE {\bf take} $g\in T$ s.t. $\lm(h)$ is $\cJ$-reducible modulo $\pol(g)$
       \IF{$\lm(h) \neq \lm(\anc(p))$}
         \IF{$\pol(p)=\theta \circ \pol(f)$ with $f\in T$, $\theta \in DP(f,T)$}
          \IF{{\it Criteria}$(p,g)$}
            \RETURN $0$
          \ENDIF
         \ENDIF
       \ELSE
         \WHILE{$h\neq 0$ {\bf and} $\exists g\in T$ s.t. $\lm(h)$ is ${\cJ}$-reducible by $q:=\pol(g)$}
            \STATE $h:=h-\lc(h)\,\theta \circ (q/\lc(q))$ with $\theta\in {\cJ}(q,T)$ and $\lm(h)=\theta \circ \lm(g)$
         \ENDWHILE
       \ENDIF
    \ENDIF \ \ {\bf return $h$}
\end{algorithmic}
\end{algorithm}

If all the criteria are false, then the ${\cJ}$-head reduction of
$h$ is done by the {\bf while-}loop 14--16 in accordance with the
def\/inition of the head reduction in~(\ref{elem_J_head_red}).

The last algorithm {\bf TailNormalForm} completes $\cJ$-reduction
of the ${\cJ}$-head reduced polynomial in the input triple by
performing its $\cJ$-tail reduction. This algorithm is invoked at
step~13 of the main algorithm {\bf Gr\"{o}bnerBasis}. The tail
reduction is performed in the {\bf while-}loop as a~sequence of
elementary reductions~(\ref{elem_J_red}).

\begin{algorithm}{TailNormalForm$(p,T,\prec)$}
\begin{algorithmic}[1]
\INPUT $p$, a triple such that
   $\pol(p)=HNF_{\cJ}(p,T)$; $T$, a set of triples; $\prec$, a ranking
\OUTPUT $h=NF_{\cJ}(p,T)$, the (full) ${\cJ}$-normal form of
$\pol(p)$ modulo $T$
    \STATE $G:=\{\,\pol(g)\mid g\in T\,\}$
    \STATE $h:=\pol(p)$
      \WHILE{$h$ has a term $t = a \vartheta \circ y^j$ which is ${\cJ}$-reducible modulo $G$}
            \STATE {\bf take} $g\in G$\ \ s.t.\ \ $\vartheta \circ y^j = \theta \circ lm(g)$
            \STATE $h:=h - \lc(h)\,\vartheta \circ (g/\lc(g)$
      \ENDWHILE
  \RETURN $h$
\end{algorithmic}
\end{algorithm}

\vskip 0.2cm

Because of the lack of an appropriate collection of benchmarks for
linear f\/inite-dif\/ference polynomial systems, the algorithmic
ef\/f\/iciency of algorithm {\bf Gr\"{o}bnerBasis} can be
indirectly analyzed by running its polynomial (non-dif\/ference)
counterpart~\cite{Gerdt'05,GB'05} for the extensive benchmarks
collection in~\cite{BM'96,JV}. Some timings for our polynomial
implementation can be found on the Web page~\cite{ginv}.

Recently, the algorithm in its dif\/ference version was
implemented in Maple~\cite{GR'05}. Just this implementation was
used for generation of linear f\/inite-dif\/ference schemes as
described in the next sections. Though one needs special and
intensive benchmarking for linear dif\/ference systems, our
f\/irst experimenting with the Maple implementation and with that
for commutative polynomials gives us a good reason to expect that
the following merits revealed for the pure polynomial
version~\cite{Gerdt'05} hold also for the dif\/ference one:
\begin{itemize}
\itemsep=0pt \item automatic avoidance of some useless reductions;
\item weakened role of the criteria: even without applying any
criteria the algorithm is reaso\-nab\-ly fast; \item smoothed
growth of intermediate coef\/f\/icients; \item fast search of a
reductor which provides the elementary Janet-like
reduction~(\ref{elem_J_red}) of a~given term. It should be noted
that there can be at most one reductor~\cite{Gerdt'05}; \item
natural and ef\/fective parallelism.
\end{itemize}

\section{Illustrative examples of PDEs}

\subsection{Laplace equation}

In this section we illustrate the approach of Section~2 to the
automatic generation of dif\/ference schemes by simplest elliptic,
parabolic and hyperbolic equations. To compute  Gr\"obner bases
providing the elimination of the partial derivatives to construct
dif\/ference schemes we used the Maple package~\cite{GR'05}
implementing the algorithms described in the previous section.

We start with the Laplace
equation~\cite{GV'96a,QV'97,Th1,Th2,Sam'01}
\begin{gather}
u_{xx}+u_{yy}=0 \label{laplace_DE}
\end{gather}
and rewrite it as the conservation law~(\ref{cons_law})
\begin{gather}
 \oint_{\Gamma} - u_y dx + u_x dy = 0. \label{laplace_IF}
\end{gather}

Now we add the relations~(\ref{rel}) for the partial derivatives
$u_x$ and $u_y$
\begin{gather}
 \int _{x_j}^{x_{j+2}} u_x dx = u(x_{j+2}, y) -
 u(x_{j}, y), \qquad
 \int _{y_k}^{y_{k+2}}  u_y dy = u(x, y_{k+2}) -
 u(x, y_{k}).
\label{laplace_rel}
\end{gather}

Thus, we obtain the system of three integral
relations~(\ref{laplace_IF}), (\ref{laplace_rel}) for three
functions
\[
u(x,y),\quad u_x(x,y),\quad u_y(x,y).
\]
To discretize this system we choose the rectangular contour of
Fig.~\ref{int_cont} on the orthogonal and uniform
grid~(\ref{discret}) with
\begin{gather}
h_1=h_2=h \label{grid_Laplace}
\end{gather}
and use the midpoint integration method for both
(\ref{laplace_IF}) and (\ref{laplace_rel}). This yields the
system:
\begin{gather*}
  -((u_y)_{j+1 \, k} - (u_y)_{j+1 \, k+2})  +
  ((u_x)_{j+2 \, k+1} - (u_y)_{j \, k+1})   = 0,
   \\
 (u_x)_{j+1 \, k} \cdot 2\,  h = u_{j+2 \, k} - u_{j \, k},
    \\
 (u_y)_{j \, k+1} \cdot 2\,  h = u_{j \, k+2} - u_{j \, k}.
\end{gather*}
Rewritten in terms of dif\/ference polynomials in the ring
$\Q\{u,u_x,u_y\}$ (see Section~2) it reads:
\begin{gather*}
    (\theta_x \theta_y^2 - \theta_x) \circ u_y  + (\theta_x^2 \theta_y - \theta_y) \circ u_x = 0,
\\
    2\,  h\, \theta_x \circ u_x   - (\theta_x^2 - 1) \circ u = 0 ,   \\
    2\,  h\, \theta_y \circ u_y  - (\theta_y^2 - 1) \circ u = 0 .
\end{gather*}
Computation of the Gr\"obner basis for the elimination ranking
(Section~2) with $u_x \succ u_y \succ u $ and $\theta_x \succ
\theta_y$ gives:
\begin{gather*}
   \theta_x \circ u_x - \frac{1}{2\,  h}\,(\theta_x^2 - 1) \circ u = 0, \\
   \theta_y \circ u_x + \theta_x \circ u_y -
   \frac{1}{2\,  h}\, (\theta_x\theta_y((\theta_x^2 - 1)+(\theta_y^2 - 1))) \circ u = 0,   \\
   \theta_x^2 \circ u_y -
   \frac{1}{2\,  h}\, (\theta_x^2\theta_y((\theta_x^2 - 1)+(\theta_y^2 - 1)) - \theta_y(\theta_x^2 - 1))
    \circ u = 0 ,  \\
    \theta_y \circ u_y - \frac{1}{2\,  h}\, (\theta_y^2 - 1) \circ u = 0 , \\
   \frac{1}{2\,  h}\,(\theta_x^4 \theta_y^2 + \theta_x^2 \theta_y^4 - 4\theta_x^2 \theta_y^2 +\theta_x^2 + \theta_y^2)
    \circ u = 0.
\end{gather*}

The latter equation with eliminated $u_x$ and $u_y$ is the
standard dif\/ference scheme with the central approximation of the
second-order derivatives in~(\ref{laplace_DE}) written in double
nodes
\[
 \frac{u_{j+2 \, k} - 2\, u_{j \, k} + u_{j - 2 \, k}}{4
h^{\,2}}
 +\frac{u_{j \, k+2} - 2\, u_{j \, k} + u_{j \, k-2}}{4
h^{\,2}}=0.
\]

Similarly, the trapezoidal integration rule for
relations~(\ref{laplace_rel}) generates the same dif\/ference
scheme but written in ordinary nodes
\begin{gather*}
 \frac{u_{j+1 \, k} - 2\, u_{j \, k} + u_{j - 1 \, k}}{
h^{\,2}}
 +\frac{u_{j \, k+1} - 2\, u_{j \, k} + u_{j \, k-1}}{
h^{\,2}}=0.
\end{gather*}

\subsection{Heat equation}

Consider now the heat
equation~\cite{GV'96a,QV'97,Th1,Th2,Sam'01,MM'05}
\[
u_{t} + \alpha u_{xx}=0
\]
in its conservation law form
\begin{gather}
 \oint_{\Gamma} -  \alpha u_x dt + u dx = 0. \label{int_HE}
\end{gather}
The integrand in~(\ref{int_HE}) contains the only partial
derivative $u_x$. Hence we add the single integral relation
\begin{gather}
 \int_{x_j}^{x_{j+1}} u_x dx = u(x_{j+1}, t) - u(x_{j}, t). \label{IR_for_HE}
\end{gather}

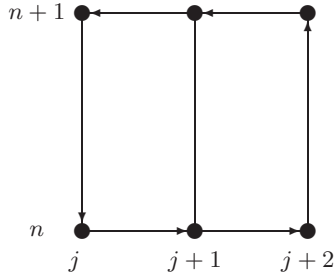
\begin{figure}[t]
\begin{center}
\unitlength=1.00mm \special{em:linewidth 0.4pt}
\linethickness{0.4pt}
\begin{picture}(110.00,34.00)(40, 4)
\put(79.00,35.00){\vector(0,-1){28.00}}
\put(79.00,35.00){\circle*{2.00}}
\put(109.00,35.00){\circle*{2.00}}
\put(109.00,6.00){\circle*{2.00}} \put(79.00,6.00){\circle*{2.00}}
\put(109.00,6.00){\vector(0,1){28.00}}
\put(73.00,6.00){\makebox(0,0)[cc]{\footnotesize $n$}}
\put(73.00,35.00){\makebox(0,0)[cc]{\footnotesize $n+1$}}
\put(78.00,2.00){\makebox(0,0)[cc]{\footnotesize $j$}}
\put(109.00,2.00){\makebox(0,0)[cc]{\footnotesize $j+2$}}
\put(94.00,6.00){\circle*{2.00}} \put(94.00,35.00){\circle*{2.00}}
\put(79.00,6.00){\vector(1,0){14.00}}
\put(94.00,6.00){\vector(1,0){14.00}}
\put(109.00,35.00){\vector(-1,0){14.00}}
\put(94.00,35.00){\vector(-1,0){14.00}}
\put(94.00,2.00){\makebox(0,0)[cc]{\footnotesize $j+1$}}
\put(94.00,6.00){\line(0,1){29.00}}
\end{picture}
\caption{Integration contour for heat
equation.}\label{HE_int_cont}
\end{center}\vspace{-5mm}
\end{figure}

Again we discretize $u(x,t)$ and $u_x(x,t)$ on the orthogonal and
uniform grid with the spatial mesh step $h$ and the temporal mesh
step $\tau$, and choose the contour shown
in~Fig.~\ref{HE_int_cont}.
Then, applying the midpoint rule for the contour integral and the
trapezoidal rule for the relation integral we f\/ind two
dif\/ference equations for two indeterminates $u, u_x$ in the form
\begin{gather*}
    \alpha \frac{\tau}{2}\, (1 + \theta_t - \theta_x^2 - \theta_t\theta_x^2) \circ u_x
    - 2\, h \,(\theta_x\theta_t - \theta_x) \circ u = 0, \\
    \frac{h}{2}\, (\theta_x + 1) \circ u_x  - (\theta_x - 1) \circ u = 0.
\end{gather*}

By elimination of $u_x$ by means of the Gr\"obner basis with
$u_x\succ u$ we obtain the famous Crank--Nicholson
scheme~\cite{GV'96a,QV'97,Th1,Th2,Sam'01,MM'05,GV'96b}
\begin{gather*}
 \frac{u_{j}^{n+1} - u_{j}^{n}}{ \tau}
 +  \alpha
\frac{(u_{j+1}^{n+1} - 2\, u_{j}^{n+1} + u_{j-1}^{n+1})
  +(u_{j+1}^{n} - 2\, u_{j}^{n} + u_{j-1}^{n})}{2\, h^{\,2}}=0.
\end{gather*}
The same scheme is obtained for the midpoint integration method
applied to~(\ref{IR_for_HE}).

\subsection{Wave equation}

The wave equation~\cite{GV'96a,QV'97,Th1,Th2,Sam'01,MM'05}
\[
u_{tt} - u_{xx}=0
\]
in the conservation law form is given by
\[
 \oint_{\Gamma}  u_x dt + u_t dx = 0.
\]
Choosing the same grid with the mesh steps~(\ref{grid_Laplace}),
the contour of Fig.~\ref{int_cont} and integral
relations~(\ref{laplace_rel}) as are used in Section~5.2 for the
Laplace equation~(\ref{laplace_IF}) and applying the midpoint rule
for the contour integral and the trapezoidal rule for the integral
relations we obtain the operator equations
\begin{gather*}
    (\theta_x - \theta_x \theta_t^2) \circ u_t  + (\theta_x^2 \theta_t - \theta_t) \circ u_x = 0,
\\
    \frac{h}{2}\,(\theta_x + 1) \circ u_x  - (\theta_x - 1)
\circ u =0, \\
    \frac{h}{2}\, (\theta_t + 1) \circ u_t   - (\theta_t - 1)
\circ u = 0.
\end{gather*}
The Gr\"obner basis method yields the standard dif\/ference scheme
\begin{gather*}
 u^{n+1}_{j} + u^{n-1}_{j}- u^{n}_{j+1} - u^{n}_{j-1}=0.
\end{gather*}

\subsection{Advection equation}

Consider now a simple form of the Advection (or convection or
one-way wave)
equation~\cite{GV'96a,QV'97,Th1,Th2,Sam'01,MM'05,GV'96b}
\begin{gather}
u_{t} + \nu\, u_{x}=0,\qquad \nu={\rm const}. \label{AE}
\end{gather}

Being of f\/irst order, the equation~(\ref{AE}) has already the
conservation law form~(\ref{cons_law}). By this reason, to
generate a dif\/ference scheme we shall not convert the equation
into the integral form~(\ref{int_cons_law}). In the latter case
one has nothing to eliminate. Instead, we consider
equation~(\ref{AE}) together with the integral relations:
\begin{gather}
    u_{t} + \nu\, u_{x}=0, \nonumber\\
    \int_{t_1}^{t_2}u_t dt=u(t_2,x)-u(t_1,x),\qquad
    \int_{x_1}^{x_2}u_x dx=u(t,x_2)-u(t,x_1).
 \label{AE_system}
\end{gather}
Discretization of $u$, $u_t$ and $u_x$ on the orthogonal and
uniform grid with the mesh steps $h$ and~$\tau$ and the explicit
integration formula for the upper integral relation
in~(\ref{AE_system}) together with the midpoint integration rule
for the lower relation give the dif\/ference system:
\begin{gather}
    u_{t} + \nu\, u_{x}=0, \qquad
    \tau \, u_t  = (\theta_t-1) \circ u,\qquad
    2\, h\,\theta_x \circ u_x = (\theta_x^2-1) \circ u.
 \label{disc_AE}
\end{gather}

Let us apply the operator $\theta_x$ to both sides of the middle
equation in~(\ref{disc_AE}) and then use the Lax method, that is,
replace $\theta_x$ with $(\theta_x^2+1)/2$ in the second term of
the right-hand side. This replacement yields
\begin{gather*}
    u_{t} + \nu\, u_{x}=0, \\
    \theta_x \circ u_t \cdot \tau - \left(\theta_t\theta_x-\frac{\theta_x^2+1}{2}\right) \circ u=0,\\
    2\,\theta_x \circ u_x \cdot h-(\theta_x^2-1) \circ u=0.
\end{gather*}
The lexicographical Gr\"obner basis for the elimination ranking
$u_t \succ u_x \succ u $ with $\theta_t \succ \theta_x$, is
\begin{gather*}
    u_{t} + \nu\, u_{x}=0, \\
    2h\,\theta_x \circ u_x  - (\theta_x^2-1) \circ u=0,\\
    (2h\,\theta_x \theta_t   - h\,(\theta_x^2+1) + \tau (\theta_x^2-1)\cdot \nu
    ) \circ  u=0.
\end{gather*}
Its last element gives the scheme:
\[
 u^{n+1}_{j+1} =\frac{u^{n}_{j+2}+u^n_j}{2}-
 \frac{\nu \tau(u^n_{j+2}-u^n_j)}{2\,h}.
\]

\section{Burgers equation}

\subsection{Conservation law form}

Consider Burgers equation~\cite{QV'97,Th1,MM'05,GV'96b} in the
form
\begin{gather}
u_t + f_x = \nu\ u_{xx}, \qquad \nu = \mathrm{const}, \label{BE}
\end{gather}
where we replaced $u^2$ by the f\/lux function $f$ in order to
avoid computation of nonlinear dif\/ference Gr\"obner bases. $\nu$
is called the viscosity. This equation exhibits some dif\/f\/icult
features from the point of view of simple f\/inite dif\/ference
schemes due to the term $f=u^2$. Let us, f\/irst, convert
equation~(\ref{BE}) into the conservation law form
\[
 \oint_{\Gamma}  (\nu u_x - f) dt + u dx = 0.
\]
Then choose the contour of Fig.~\ref{int_cont} and add the
integral relation
\[
 \int_{x_j}^{x_{j+2}}  u_x dx = u(x_{j+2}, t) -
 u(x_{j}, t).
\]

Denoting as above the temporal and spatial mesh steps by $\tau$
and $h$, and applying the midpoint integration rule we obtain the
system:
\begin{gather*}
     h \,(\theta_x \theta_t^2 - \theta_x) \circ u  - \tau \,(\theta_x^2 \theta_t - \theta_t)\circ (\nu u_x - f)
     = 0, \\
    2\, h\, \theta_x \circ u_x   - (\theta_x^2 - 1) \circ u = 0 .
\end{gather*}
Its Gr\"obner basis form for the elimination order with $u_x \succ
u \succ f$ and $\theta_t \succ \theta_x$ is given by
\begin{gather*}
   2\,\nu \tau h\, \theta_t\circ u_x + 2\, h^2\theta_x(\theta^2_t-1)\circ u +2\,\tau h
   \theta_t(\theta_x^2-1)\circ f -
    \nu \tau \theta_t\theta_x(\theta_x^2-1)\circ u=0, \\
      2\,  h\, \theta_x \circ u_x   - (\theta_x^2 - 1) \circ u = 0,\\
    2\, h^2 \theta_x^2(\theta_t^2 - 1) \circ u  - \nu \tau\theta_t (\theta_x^4 - 2\, \theta_x^2 +1)\circ u +
    2\, \tau  h \,\theta_t\theta_x (\theta_x^2-1)\circ f = 0.
\end{gather*}
The obtained dif\/ference scheme
\begin{gather}
\frac{u^{n+2}_{j+2}-u^n_{j+2}}{\tau}-\nu
\frac{u^{n+1}_{j+4}-2u^{n+1}_{j+2}+u^{n+1}_j}{2h^{\,2}}+
 \frac{f^{n+1}_{j+3}-f^{n+1}_{j+1}}{h}=0. \label{BE_FTFS}
\end{gather}
is the standard explicit scheme with forward time and forward
space dif\/ferencing. It is well-known that schemes of this type
are unstable~\cite{Th1,MM'05,GV'96b}. Furthermore, by using
implicit schemes one can provide the von Neumann
stability\footnote{For example, if one uses the central
dif\/ferencing in the last term of~(\ref{BE_FTFS}), and the
Crank--Nicolson approach to the second (dif\/fusion) term.}.
However, all such schemes are usually not very satisfactory when
one considers non-smooth or discontinuous solutions (shock waves)
of Burgers equation.

\subsection{Lax method}
To exploit more f\/lexibility and freedom in our dif\/ference
elimination approach to generation of f\/inite-dif\/ference
schemes, we go back to the original dif\/ferential
equation~(\ref{BE}) and consider it together with the integral
relations providing the elimination. For discretization of the
relations we combine the midpoint rule for integration over $x$
with the explicit integration over $t$ and apply the Lax method to
the last integration:
\begin{alignat}{3}
& u_t+f_x=\nu u_{xx}, && (u_t)^n_j+(f_x)^n_j=\nu (u_{xx})^n_j,&\nonumber\\
& \int u_t dt = u , && u_t \tau = u_j^{n+1}-\frac{u^n_{j+2}+u^n_j}{2},& \nonumber\\
&     \int f_x dx = f, &\qquad \Longrightarrow\qquad &  2\,h (f_x)^n_{j+1} = f_{j+2}^n-f_j^n , & \label{bya:1}\\
&  \int u_x dx = u , && 2\,h (u_x)^n_{j+1} = u_{j+2}^n-u_j^n,& \nonumber\\
&  \int u_{xx} dx = u_x, && 2\,h (u_{xx})^n_{j+1} =
(u_x)_{j+2}^n-(u_x)_j^n.& \nonumber
\end{alignat}
The Gr\"obner basis based elimination with $u_{xx}\succ u_t\succ
u_x\succ f_x \succ u \succ f$ yields the scheme
\begin{gather}
\frac{2\,u^{n+1}_{j+2} - (u^{n}_{j+3} + u^{n}_{j+1})}{2\tau}
    + \frac{f^{n}_{j+3} - f^{n}_{j+1}}{2h} - \nu\frac{u^{n}_{j+4}-2\,u^{n}_{j+2}+u^{n}_{j}}{4h^{\,2}}=0.
\label{bya:2}
\end{gather}

One can also use the trapezoidal rule for the spatial
integrations. This derives other schemes. Since there are three
spatial integrals in~(\ref{bya:1}), by selecting either the
midpoint or the trapezoidal rule for these integrals, we obtain
eight possible variants of the dif\/ference schemes. Our
computation with the Gr\"obner bases reveals seven dif\/ferent
schemes. Apart from~(\ref{bya:2}) there are
\begin{gather}
  \frac{2(u^{n+1}_{j+2}+u^{n+1}_{j+1}) - (u^{n}_{j+3}+u^{n}_{j+2}+u^{n}_{j+1}+u^{n}_{j})}{4\tau}
    + \frac{(f^{n}_{j+3}+f^{n}_{j+2}) - (f^{n}_{j+1}+f^{n}_{j})}{4h}\nonumber\\
\qquad{}   =
\nu\frac{(u^{n}_{j+3}-u^{n}_{j+2})-(u^{n}_{j+1}-u^{n}_{j})}{2h^2},
  \label{bya:112}
\\
  \frac{2u^{n+1}_{j+1} - (u^{n}_{j+2} + u^{n}_{j})}{2\tau}
    + \frac{f^{n}_{j+2} - f^{n}_{j}}{2h}
    = \nu\frac{u^{n}_{j+2}-2u^{n}_{j+1}+u^{n}_{j}}{h^2},
  \label{bya:121}
\\
  \frac{2(u^{n+1}_{j+3}+2u^{n+1}_{j+2}+u^{n+1}_{j+1}) - (u^{n}_{j+4}+2u^{n}_{j+3}+2u^{n}_{j+2}+2u^{n}_{j+1}+u^{n}_{j})}{8\tau}
   \nonumber\\
 \qquad {} + \frac{(f^{n}_{j+4}+2f^{n}_{j+1}) - (2f^{n}_{j+1}+f^{n}_{j})}{8h}
    = \nu\frac{u^{n}_{j+3}-2u^{n}_{j+2}+u^{n}_{j+1}}{h^2},
  \label{bya:122}
\\
  \frac{2(u^{n+1}_{j+3}+u^{n+1}_{j+2}) - (u^{n}_{j+4}+u^{n}_{j+3}+u^{n}_{j+2}+u^{n}_{j+1})}{4\tau}
    + \frac{f^{n}_{j+3} - f^{n}_{j+2}}{h}\nonumber\\
\qquad{}   =
\nu\frac{((u^{n}_{j+5}+u^{n}_{j+4})-2u^{n}_{j+3})-(2u^{n}_{j+2}-(u^{n}_{j+1}+u^{n}_{j}))}{8h^2},
  \label{bya:211}
\\
  \frac{2(u^{n+1}_{j+2}+u^{n+1}_{j+1}) - (u^{n}_{j+3}+u^{n}_{j+2}+u^{n}_{j+1}+u^{n}_{j})}{4\tau}
    + \frac{f^{n}_{j+2} - f^{n}_{j+1}}{h}\nonumber\\
\qquad{}
    = \nu\frac{(u^{n}_{j+3}-u^{n}_{j+2})-(u^{n}_{j+1}-u^{n}_{j})}{2h^2},
  \label{bya:212:221}\\
  \frac{2(u^{n+1}_{j+3}+2u^{n+1}_{j+2}+u^{n+1}_{j+1})
  - (u^{n}_{j+4}+2u^{n}_{j+3}+2u^{n}_{j+2}+2u^{n}_{j+1}+u^{n}_{j})}{8\tau}
    + \frac{f^{n}_{j+3} - f^{n}_{j+1}}{2h}\nonumber\\
\qquad{}    = \nu\frac{u^{n}_{j+3}-2u^{n}_{j+2}+u^{n}_{j+1}}{h^2}.
  \label{bya:222}
\end{gather}
Just the scheme~(\ref{bya:212:221}) is obtained twice in the
course of generating eight schemes.

\subsection[Two-step Lax-Wendroff method]{Two-step Lax--Wendrof\/f method}
Our Gr\"obner basis based technique can also be applied to
generate other types of dif\/ference schemes. For example, one can
generate two-step Lax--Wendrof\/f schemes~\cite{RM'67}. Let
$\overline{u}$ and $\overline{f}$ denote the values of $u$ and $f$
on the intermediate time levels. Then, applying again the midpoint
rule for the spatial integrals, gives the following dif\/ference
system:
\begin{gather*}
    {u_t\,}^{n}_{j} + {f_x\,}^{n}_{j} = \nu\ {u_{xx}\,}^{n}_{j}, \\
    {u_t\,}^{n}_{j}\,\tau = {\overline{u}\,}^{n+1}_{j} - \frac{u^{n}_{j+2} + u^{n}_{j}}{2}, \\
    2{f_x\,}^{n}_{j+1}\,h = f^{n}_{j+2} - f^{n}_{j}, \\
    2{u_x\,}^{n}_{j+1}\,h = u^{n}_{j+2} - u^{n}_{j}, \\
    2{u_{xx}\,}^{n}_{j+1}\,h = {u_x\,}^{n}_{j+2} - {u_x\,}^{n}_{j}, \\
    {\overline{u}_t\,}^{n}_{j} + {\overline{f}_x\,}^{n}_{j} = \nu\ {\overline{u}_{xx}\,}^{n}_{j}, \\
    {\overline{u}_t\,}^{n}_{j}\,\tau = u^{n+1}_{j} - u^{n}_{j}, \\
    2{\overline{f}_x\,}^{n}_{j+1}\,h = \overline{f}^{n}_{j+2} - \overline{f}^{n}_{j},\\
    2{\overline{u}_x\,}^{n}_{j+1}\,h = \overline{u}^{n}_{j+2} - \overline{u}^{n}_{j}, \\
    2{\overline{u}_xx\,}^{n}_{j+1}\,h = {\overline{u}_x\,}^{n}_{j+2} - {\overline{u}_x\,}^{n}_{j}.
\end{gather*}
For the elimination ranking with
\begin{displaymath}
u_{xx} \succ \overline{u}_{xx} \succ u_x \succ \overline{u}_{x}
\succ u_t \succ \overline{u}_{t} \succ f_x \succ \overline{f}_{x}
\succ f \succ u \succ \overline{f} \succ \overline{u}
\end{displaymath}
the Gr\"obner basis contains the Lax--Wendrof\/f scheme
\begin{gather}
  \frac{\overline{u}^{n+1}_{j+2} - (u^{n}_{j+3} + u^{n}_{j+1})}{2\,\tau}
    + \frac{f^{n}_{j+3} - f^{n}_{j+1}}{2\,h}
    = \nu\frac{u^{n}_{j+4}-2u^{n}_{j+2}+u^{n}_{j}}{4h^2}, \nonumber\\
  \frac{u^{n+1}_{j+3} - \overline{u}^{n}_{j+2}}{2\,\tau}
    + \frac{\overline{f}^{n}_{j+3} - \overline{f}^{n}_{j+1}}{2\,h}
    = \nu\frac{\overline{u}^{n}_{j+4}-\overline{u}^{n}_{j+2}+\overline{u}^{n}_{j}}{4\,h^2}.
\label{bya:17}
\end{gather}
With all possible combinations of the trapezoidal and midpoint
rules one obtains 49 dif\/ferent Lax--Wendrof\/f schemes.

\subsection[Differential approximation]{Dif\/ferential approximation}
To analyze properties of a dif\/ference scheme it can be useful to
compute its dif\/ferential approximation~\cite{SY'85} that is
often called the modif\/ied equation(s) of the dif\/ference
scheme. There are whole classes of dif\/ferent schemes for which
their stability properties can be obtained with the aid of the
dif\/ferential approximation~\cite{Str'89}. For all that, in many
cases, the computation can be easily done with modern computer
algebra software. In our research we use Maple~\cite{Maple}.
Consider, for example, the schemes~(\ref{bya:2})--(\ref{bya:222}).

Their dif\/ferential approximation for $f=u^2$ and with collection
of the coef\/f\/icients at $\tau$, $h^2$, $h^2/\tau$ is given by:
\begin{gather*}
u_t + u_xu - \nu\, u_{xx}=\left(-\frac
12\,\nu^2u_{xxxx}+(u_{xxx}u+2\,u_{xx}u_{x})\nu-u_{x}^2u
-\frac 12 \,u^2u_{xx}\right)\tau\\
\phantom{u_t + u_xu - \nu\, u_{xx}=}+(*)h^2 + \frac 12\,
u_{xx}\frac{h^2}{\tau} + \cdots.
\end{gather*}

The schemes~(\ref{bya:2})--(\ref{bya:222}) dif\/fer in the
coef\/f\/icient $(*)$ at $h^2$ only. These coef\/f\/icients are as
follows
\begin{alignat*}{3}
& (32) &\qquad & -\frac 16 \,u_{xxxx}\nu+\frac 13 \,u_{xxx}u+\frac 12 \,u_{xx}u_{x},& \\
& (33) && -\frac 16 \,u_{xxxx}\nu+\frac{1}{12}\,u_{xxx}u-\frac{1}{4}\,u_{xx}u_{x}, &\\
& (34) && -\frac{5}{12}\,u_{xxxx}\nu+\frac{1}{3}\,u_{xxx}u+\frac{1}{2}\,u_{xx}u_{x}, & \\
& (35) && -\frac{1}{6}\,u_{xxxx}\nu-\frac{1}{6}\,u_{xxx}u-u_{xx}u_{x}, & \\
& (36) &&  \frac{1}{12}\,u_{xxxx}\nu+\frac{1}{3}\,u_{xxx}u+\frac{1}{2}\,u_{xx}u_{x}, & \\
& (37) && -\frac{1}{6}\,u_{xxxx}\nu+\frac{1}{3}\,u_{xxx}u +\frac{1}{2}\,u_{xx}u_{x}, & \\
& (38) && -\frac{1}{6}\,u_{xxxx}\nu -
\frac{1}{4}\,u_{xx}u_{x}+\frac{1}{12}\,u_{xxx}u . &
\end{alignat*}

Thereby, comparison of dif\/ferential approximations for
schemes~(\ref{bya:2})--(\ref{bya:222}) shows that
\begin{itemize}\itemsep=0pt
\item all the schemes provide the same order of approximation in
$\tau$, $h$; \item they have identical linear numerical
dissipation (viscosity)~\cite{QV'97,GV'96b} determined by
$u_{xx}h^2/(2\,\tau)$; \item the schemes possess similar
dispersion properties with distinction in the rational
coef\/f\/icients of the dif\/ferential polynomial in $u$
multiplied by $h^2$.
\end{itemize}

As to scheme~(\ref{BE_FTFS}), the right-hand side of its
dif\/ferential approximation reads
\begin{gather*}
\left(-\frac{1}{2}\nu^2u_{xxxx}+(u_{xxx}u+2u_{xx}u_{x})\nu-u_{x}^2u-\frac{1}{2} u^2u_{xx}\right)\tau\\
\qquad{}+ \left(\frac{1}{3} u_{xxxx}\nu-\frac{1}{6}
u_{xxx}u-\frac{1}{2} u_{xx}u_{x}\right)h^2 + \cdots.
\end{gather*}
This explicitly shows instability of the scheme which does not
yield linear numerical viscosity.

We obtained also analogous results on stability and on close
properties for the dif\/ferent Lax--Wendrof\/f schemes of
type~(\ref{bya:17}) and its variations due to the choice of
dif\/ferent numerical integration rules for the spatial integrals.

\subsection{Godunov method}
It is especially dif\/f\/icult to simulate numerically nonsmooth
and discontinuous solutions which are among most interesting
problems in computational f\/luid and gas
dynamics~\cite{GR'87,QV'97,Th1,MM'05}. Most of the known
dif\/ference schemes fail to handle these singularities. The most
appropriate numerical approach to such problems was developed by
Godunov~\cite{GR'87,God'59} and based on solving a local Riemann
problem~\cite{QV'97,Th2} as a cornerstone of the Godunov's scheme
generation. There are special numerical Riemann solvers, for
example~\cite{Toro'97}, designed for these purposes and for
application to computational f\/luid dynamics.

Instead of the use of numerical Riemann solvers,  we apply the
Gr\"obner bases technique to generate the Godunov-type scheme for
inviscid Burgers equation when $\nu=0$ in~(\ref{BE}). For this
purpose we discretize the corresponding system in (\ref{bya:1}) in
the following way
\begin{gather}
    {u_t\,}^{n}_{j} + {f_x\,}^{n}_{j} = 0, \nonumber\\
    {u_t\,}^{n}_{j}\,\tau = u^{n+1}_{j} - u^{n}_{j},  \nonumber\\
    ({f_x\,}^{n}_{j}\,h - (f^{n}_{j+1} - f^{n}_{j}))
    ({f_x\,}^{n}_{j+1}\,h - (f^{n}_{j+1} - f^{n}_{j}))=0, \nonumber\\
    2\,{u_x\,}^{n}_{j+1}\,h = u^{n}_{j+2} - u^{n}_{j}, \nonumber\\
    2\,{u_{xx}\,}^{n}_{j+1}\,h = {u_x\,}^{n}_{j+2} - {u_x\,}^{n}_{j}.\label{bya:25}
\end{gather}
Here, the third equation contains in its left-hand side the
product of two dif\/ferent solutions for the f\/lux function $f$
of the local Riemann problem~\cite{Toro'97}. Therefore, we add to
the system composed of the original dif\/ferential equation and
discrete forms of the integral relations for partial derivatives
$u_t$, $u_{x}$, $u_{xx}$ the nonlinear dif\/ference equation on
$f$ and $f_x$ containing solutions of the local Riemann problem.

Since the Riemann condition on the f\/lux is now a constituent of
the dif\/ference system, an elimination of all the partial
derivatives of $u$ and $f$ gives the dif\/ference scheme
consistent with that condition. To do the elimination from the
nonlinear system~(\ref{bya:25}) we apply the Gr\"obner factoring
approach~\cite{Cz'89}: if a Gr\"obner basis contains a polynomial
which factors, then the computation is split into the computation
of two or more Gr\"obner bases corresponding to the factors. In
doing so, we choose the elimination ranking
\[
u_{xx}\succ u_x \succ u_t\succ f_x\succ f\succ u
\]
and compute two Gr\"obner bases, for every factor
in~(\ref{bya:25}). Then we compose the product of two obtained
dif\/ference polynomials in $u$ and $f$ that gives us the
Godunov-type dif\/ference scheme:
\begin{gather}
  \left(\frac{u^{n+1}_{j+2}-u^{n}_{j+2}}{\tau}
        +\frac{f^{n}_{j+2}-f^{n}_{j+1}}{h}
  \right) 
   \cdot \left(\frac{u^{n+1}_{j+2}-u^{n}_{j+2}}{\tau}
        +\frac{f^{n}_{j+3}-f^{n}_{j+2}}{h}
   \right)=0.
\label{bya:26}
\end{gather}

Below (Section~7) we compare schemes~(\ref{bya:2}), (\ref{bya:17})
and (\ref{bya:26}) by numerical simulation of a~discontinuous
solution.

\section[Falkowich-Karman equation]{Falkowich--Karman equation}

\begin{figure}[htbp]
\begin{center}
\unitlength=1.00mm \special{em:linewidth 0.4pt}
\linethickness{0.4pt}
\begin{picture}(108.00,47.00)
\put(52.00,5.00){\makebox(0,0)[cc]{\footnotesize $j+2$}}
\put(40.00,9.00){\makebox(0,0)[cc]{\footnotesize $j+1$}}
\put(28.00,13.00){\makebox(0,0)[cc]{\footnotesize $j$}}
\put(65.00,5.00){\makebox(0,0)[cc]{\footnotesize $k$}}
\put(74.00,8.00){\makebox(0,0)[cc]{\footnotesize $k+1$}}
\put(82.00,12.00){\makebox(0,0)[cc]{\footnotesize $k+2$}}
\put(28.00,18.00){\line(3,-1){31.00}}
\put(59.00,7.67){\line(5,2){22.00}}
\put(81.00,16.33){\line(-3,1){32.00}}
\put(49.00,27.00){\line(-5,-2){21.00}}
\put(41.00,23.00){\line(3,-1){30.00}}
\put(28.00,38.00){\line(3,-1){31.00}}
\put(59.00,27.67){\line(5,2){22.00}}
\put(81.00,36.33){\line(-3,1){32.00}}
\put(49.00,47.00){\line(-5,-2){21.00}}
\put(41.00,43.00){\line(3,-1){30.00}}
\put(44.00,13.00){\circle*{2.00}}
\put(65.00,21.00){\circle*{2.00}}
\put(41.00,23.00){\circle*{2.00}}
\put(72.00,13.00){\circle*{2.00}}
\put(41.00,43.00){\circle*{2.00}}
\put(56.00,18.00){\circle*{2.00}}
\put(72.00,33.00){\circle*{2.00}}
\put(44.00,13.00){\line(5,2){21.00}}
\put(72.00,13.00){\line(0,1){19.00}}
\put(41.00,23.00){\line(0,1){20.00}}
\put(22.00,20.00){\makebox(0,0)[cc]{\footnotesize $n$}}
\put(22.00,40.00){\makebox(0,0)[cc]{\footnotesize $n+1$}}
\end{picture}
\vspace{-5mm} \caption{Cell for Falkowich--Karman equation.}
   \label{int_volume}
\end{center}
\vspace{-7mm}
\end{figure}
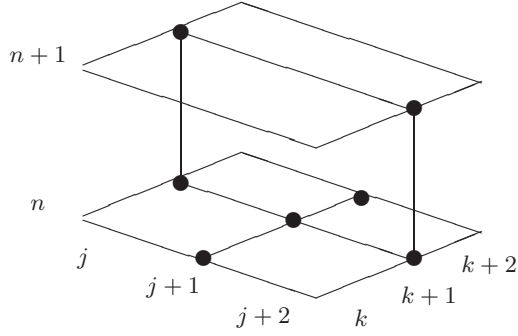

Consider now the nonlinear two-dimensional Falkowich--Karman
equation~\cite{FKar}
\begin{gather}
\varphi_{xx}(K - (\gamma +1)\varphi_x)+\varphi_{yy}=0
\label{FKE_ST}
\end{gather}
describing transonic f\/low in gas dynamics in its non-stationary
form
\begin{gather}
\varphi_{xx}(K - (\gamma
+1)\varphi_x)+\varphi_{yy}-2\varphi_{xt}-\varphi_{tt}=0.
\label{FKE}
\end{gather}
This form can be used to f\/ind a stationary solution by the
steady-state method. We rewrite equation~(\ref{FKE}) into the
conservation law form
\[
 \int_{t_{n}}^{t_{n+1}}\left(
  \oint_{\Gamma} - \varphi_{y} dx +
   \varphi_{x}\left(K - \frac{(\gamma +1)}{2} \varphi_x\right) dy \right) dt -
 \int_{x_{j}}^{x_{j+2}} \int_{y_{k}}^{y_{k+2}}
 (2\varphi_{x} + \varphi_{t})\Big|_{t_{n}}^{t_{n+1}} dx dy
   = 0.
\]
Here we use a grid decomposition of the three-dimensional domain
in $(x,y,t)$ into elementary volumes. Fig.~\ref{int_volume} shows
an elementary volume.

Again we add the integral relations for partial derivatives with
the use of the trapezoidal integration rule for $\varphi_{x}$,
$\varphi_{y}$ and the midpoint rule for $\varphi_{t}$.

Then we obtain the nonlinear operator equations:
\begin{gather}
    \left(-(\theta_x - \theta_x \theta_y^2) \circ \varphi_y  + (\theta_x^2 \theta_y - \theta_y) \circ
    \left(\varphi_x\left(K - \frac{(\gamma +1)}{2}\, \varphi_x\right)\right)\right)\cdot 2\,
h\,\tau \nonumber\\
\qquad{}
   - (\theta_t-1)(\theta_x^2 \theta_y - \theta_y) \circ 2\, \varphi \cdot 2\, h -
    \theta_x \theta_y \circ \varphi_t \cdot 4\, h^2  = 0,\nonumber \\
    (\theta_x + 1) \circ \varphi_x \cdot \frac{ h}{2} = (\theta_x -
1) \circ \varphi, \nonumber \\
    (\theta_y + 1) \circ \varphi_y \cdot \frac{ h}{2} = (\theta_y -
1) \circ \varphi,\nonumber \\[0.1cm]
    \theta_t \circ \varphi_t \cdot 2\,\tau = (\theta_t^2 - 1) \circ \varphi.
 \label{DF_of_FKE}
\end{gather}

Because of nonlinearity in the initial dif\/ferential
equation~(\ref{FKE}), the dif\/ference system obtained is also
nonlinear. By this reason the Maple package~\cite{GR'05}
implementing algorithm {\bf Gr\"{o}bnerBasis} (Section~4) is
inapplicable to~(\ref{DF_of_FKE}). Since there is no software
available for computing dif\/ference Gr\"obner bases for nonlinear
systems, we had to perform hand computations in accordance with
the above described algorithm and with an assistance of Maple to
check some intermediate results.

In these calculations we used the lexicographical ranking such
that $\varphi_x \succ \varphi_y  \succ \varphi_t \succ \varphi $
and $\theta_x \succ \theta_y \succ \theta_t$. The resulting
Gr\"obner basis has the form:
\begin{gather*}
    (\theta_x - 1)^2\theta_y \circ \varphi \cdot (\gamma +1) \varphi_x
    \frac{\tau}{h}
    + \theta_x \theta_y \circ \varphi_t \cdot h    \\ \qquad
  {}=  ((\theta_x - 1)^2\theta_y \circ \varphi \cdot
     (K - (\theta_x - 1)^2\theta_y \circ (\gamma +1) \varphi)  +
    \theta_x(\theta_y -1)^2 \circ \varphi)\frac{\tau}{h}\\
\qquad{}   - (\theta_x^2 \theta_y - \theta_y)(\theta_t-1) \circ \varphi, \\
    (\theta_y + 1) \circ \varphi_y = (\theta_y - 1) \circ \varphi \cdot
\frac{2}{h}, \\
    \theta_t \circ \varphi_t = (\theta_t^2 - 1) \circ \varphi  \cdot
\frac{1}{2\tau},    \\
   0= \theta_x(\theta_x - 1)^2\theta_y\theta_t \circ \varphi \cdot
   \Big[\Big((\theta_x - 1)^2\theta_y\theta_t \circ \varphi \cdot
     \Big(K      - (\theta_x^3 - \theta_x^2 +\theta_x -1)\theta_y\theta_t \circ \frac{(\gamma +1)}{2}
\varphi\Big) \\
    \qquad {}+ \theta_x(\theta_y -1)^2\theta_t \circ \varphi\Big)\frac{\tau}{h}
   - (\theta_x^2 - 1)\theta_y(\theta_t^2-\theta_t) \circ \varphi -
   \theta_x\theta_y(\theta_t - 1)^2 \circ \varphi \cdot \frac{h}{2\tau}\Big]
    \\ \qquad {}+
   (\theta_x - 1)^2\theta_y\theta_t \circ \varphi \cdot
   \Big[\Big(\theta_x(\theta_x - 1)^2\theta_y\theta_t \circ \varphi \cdot
     \Big(K
     - (\theta_x^3 - \theta_x^2 +\theta_x -1)\theta_y\theta_t \circ \frac{(\gamma +1)}{2}
\varphi\Big) \\
    \qquad {}+ \theta_x^2(\theta_y -1)^2\theta_t \circ \varphi\Big)\frac{\tau}{h}
   - (\theta_x^3 - \theta_x)\theta_y(\theta_t^2-\theta_t) \circ \varphi -
   \theta_x^2\theta_y(\theta_t - 1)^2 \circ \varphi \cdot \frac{h}{2\tau}\Big].
\end{gather*}
The last element is the f\/inite-dif\/ference scheme for
equation~(\ref{FKE}):
\begin{gather}
 (\varphi_{j+1 \, k}^n - 2 \varphi_{j \, k}^n + \varphi_{j - 1 \,
k}^n) \cdot
 \Big[\Big((\varphi_{j \, k}^n - 2 \varphi_{j-1 \, k}^n + \varphi_{j - 2 \,
k}^n)
 \Big(K -  \frac{(\gamma +1)}{2h}(\varphi_{j+1 \, k}^n -\varphi_{j \, k}^n \nonumber\\
\qquad{} + \varphi_{j - 1 \, k}^n - \varphi_{j - 2 \, k}^n) +
  (\varphi_{j-1 \, k+1}^n - 2 \varphi_{j-1 \, k}^n + \varphi_{j-1
\, k-1}^n)
 \Big)\frac{\tau}{h}  \nonumber\\
\qquad{} - (\varphi_{j \, k}^{n+1} - \varphi_{j-2 \, k}^{n+1}
  - \varphi_{j \, k}^n + \varphi_{j - 2 \, k}^n) -
  (\varphi_{j-1 \, k}^{n+1} - 2 \varphi_{j-1 \, k}^n + \varphi_{j -
1 \, k}^{n-1})
  \frac{h}{2\tau} \Big]  \nonumber\\
\qquad{} +(\varphi_{j \, k}^n - 2 \varphi_{j-1 \, k}^n +
\varphi_{j - 2 \, k}^n) \cdot
 \Big[\Big((\varphi_{j+1 \, k}^n - 2 \varphi_{j \, k}^n + \varphi_{j - 1 \,
k}^n)
 \Big(K -   \frac{(\gamma +1)}{2h}(\varphi_{j+1 \, k}^n -\varphi_{j \, k}^n \nonumber\\
\qquad{}  + \varphi_{j - 1 \, k}^n - \varphi_{j - 2 \, k}^n) +
  (\varphi_{j \, k+1}^n - 2 \varphi_{j \, k}^n + \varphi_{j \, k-
1}^n)
 \Big) \frac{\tau}{h}  \nonumber\\
 \qquad{} -(\varphi_{j+1 \, k}^{n+1} - \varphi_{j-1 \, k}^{n+1}
  - \varphi_{j+1 \, k}^n + \varphi_{j - 1 \, k}^n) -
  (\varphi_{j \, k}^{n+1} - 2 \varphi_{j \, k}^n + \varphi_{j \,
k}^{n-1})
  \frac{h}{2\tau} \Big] = 0.
 \label{S_FKE}
\end{gather}
By construction, this scheme is fully conservative and does not
involve switches that is typical for computing transonic
f\/lows~\cite{J'76}.

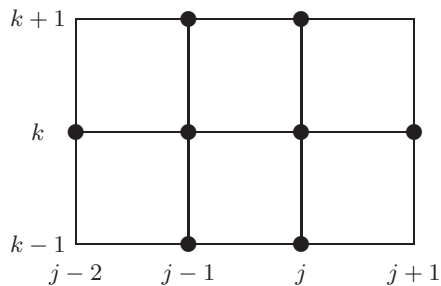
\begin{figure}[t]
\begin{center}
\unitlength=1.00mm \special{em:linewidth 0.4pt}
\linethickness{0.4pt}
\begin{picture}(120.00,38.00)(-5, 0)
\put(34.00,7.00){\line(1,0){45.00}}
\put(79.00,7.00){\line(0,1){30.00}}
\put(79.00,37.00){\line(-1,0){45.00}}
\put(34.00,37.00){\line(0,-1){30.00}}
\put(49.00,7.00){\line(0,1){30.00}}
\put(64.00,37.00){\line(0,-1){30.00}}
\put(34.00,22.00){\line(1,0){45.00}}
\put(79.00,22.00){\circle*{2.00}}
\put(64.00,22.00){\circle*{2.00}}
\put(64.00,37.00){\circle*{2.00}}
\put(49.00,37.00){\circle*{2.00}}
\put(49.00,22.00){\circle*{2.00}}
\put(34.00,22.00){\circle*{2.00}} \put(49.00,7.00){\circle*{2.00}}
\put(64.00,7.00){\circle*{2.00}}
\put(29.00,37.00){\makebox(0,0)[cc]{\footnotesize $k+1$}}
\put(29.00,22.00){\makebox(0,0)[cc]{\footnotesize $k$}}
\put(29.00,7.00){\makebox(0,0)[cc]{\footnotesize $k-1$}}
\put(34.00,3.00){\makebox(0,0)[cc]{\footnotesize $j-2$}}
\put(49.00,3.00){\makebox(0,0)[cc]{\footnotesize $j-1$}}
\put(64.00,3.00){\makebox(0,0)[cc]{\footnotesize $j$}}
\put(79.00,3.00){\makebox(0,0)[cc]{\footnotesize $j+1$}}
\end{picture}
\vspace{-3mm} \caption{Stencil for stationary Falkowich--Karman
equation.} \label{stencil}
\end{center}
\vspace{-5mm}
\end{figure}

In its stationary form scheme~(\ref{S_FKE})
\begin{gather}
 D_{xx}(\varphi_{j \, k}^n)\cdot
 \Bigg[D_{xx}(\varphi_{j-1 \, k}^n)\Bigg(K -  \frac{(\gamma +1)}{2}
   (D_{x}(\varphi_{j+1 \, k}^n) + D_{x}(\varphi_{j-1 \, k}^n))\Bigg)
   + D_{yy}(\varphi_{j-1 \, k}^n)\Bigg]\nonumber\\
\qquad{}+ D_{xx}(\varphi_{j-1 \, k}^n)\cdot
 \Bigg[D_{xx}(\varphi_{j \, k}^n)\Bigg(K -  \frac{(\gamma +1)}{2}
   (D_{x}(\varphi_{j+1 \, k}^n) + D_{x}(\varphi_{j-1 \, k}^n))\Bigg)\nonumber\\
\qquad{}   + D_{yy}(\varphi_{j \, k}^n)\Bigg]=0. \label{S_FKE_ST}
\end{gather}
is related to equation~(\ref{FKE_ST}). Here symbols $D_{x}$ and
$D_{y}$ are the forward dif\/ferencing operators and~$D_{xx}$ and
$D_{yy}$ are the central second-order dif\/ferencing operators
with respect to $x$ and $y$. The stencil for
scheme~(\ref{S_FKE_ST}) is shown in Fig.~\ref{stencil}.

It should be noted that, unlike the original dif\/ferential
equation~(\ref{FKE}) which is quadratically nonlinear, both
schemes~(\ref{S_FKE}) and (\ref{S_FKE_ST}) have the the cubic
nonlinearity in the grid function. This is the rigorous algebraic
consequence of the dif\/ference system~(\ref{DF_of_FKE}). In
accordance to the well-known fact~\cite{BAL'97}, that algebraic
elimination of variables from a nonlinear system leads generally
to increase of its degree of nonlinearity.

As an application of this scheme, in the next section we consider
an example of one-di\-men\-sio\-nal transonic f\/low with
shock-wave taken from~\cite{J'76}.

\section{Numerical experiments}

\subsection{Burgers equation}

We used schemes~(\ref{bya:2}), (\ref{bya:17}) and (\ref{bya:26})
for numerical simulation in $0<x<1$ of the following Riemann
problem for the inviscid Burgers equation~(\ref{BE}) with $\nu=0$
\begin{gather}
u_t + u_xu = 0, \label{BE_nu=0}
\end{gather}
and discontinuous initial condition $u(x,0)$:
\begin{gather}
u(x, 0)=\left\{ \begin{array}{ll}
u_l, & \  0 < x < \frac 12, \\[0.2cm]
u_r, & \   \frac 12 < x < 1.
\end{array}\right. \label{IC}
\end{gather}
In~(\ref{IC}) the initial data at $t=0$ is a piecewise-constant
function with the state $u_l$ on the left of the discontinuity
$x=0$ and the state $u_r$ on the right of the discontinuity. We
consider $\nu=0$, since in this case the problem (\ref{BE_nu=0}),
(\ref{IC}) admits the exact solution:
\begin{gather}
u(x,t)=u_lH\left(\frac{1}{2}+\frac{u_l+u_r}{2}\,t-x\right)+u_rH\left(x-\frac{1}{2}-
\frac{u_l+u_r}{2}\,t\right). \label{ES}
\end{gather}
Here $H(y)$ is the Heaviside step function~\cite{SO'87} whose
derivative is the Dirac delta function:
\begin{displaymath}
H(y)=\left\{
\begin{array}{l}
0\quad  y<0,\\[0.2cm]
\frac{1}{2}\quad  y=0,\\[0.2cm]
1\quad  y>0,
\end{array}
\right. \qquad \frac{d}{dy}H(y)=\delta (y).
\end{displaymath}
Physically, the solution~(\ref{ES}) def\/ined by the initial
condition~(\ref{IC}) represents a shock wave which moves with
constant speed $(u_l+u_r)/2$ without changing its shape.

In our numerical simulation the values of $u_l$ and $u_l$ were
chosen as 0.8 and 0.2. The  pictures below demonstrate the
numerical solution of the Riemann problem~(\ref{BE_nu=0}),
(\ref{IC}) at time $t=2/3$. Solid line shows the exact
solution~(\ref{ES}), and the numerical results are depicted by
green dots. For the ratio $\tau/h$ of mesh steps which is called
Courant (or Courant--Friedrichs--Levy) number~\cite{GR'87} we have
chosen the two values $0.9$ and $0.1$.

All schemes are numerically stable. For schemes~(\ref{bya:2})
and~(\ref{bya:17}) their stability is analytically showed by the
dif\/ferential approximation (Section~6.4). Because of the
nonlinearity in $f_x$ in the third equation of Godunov
scheme~(\ref{bya:26})  we did not compute the dif\/ferential
approximation for this scheme.

\begin{figure}[t]
\begin{minipage}[b]{7.5cm}
\centerline{\includegraphics[width=5cm,angle=270]{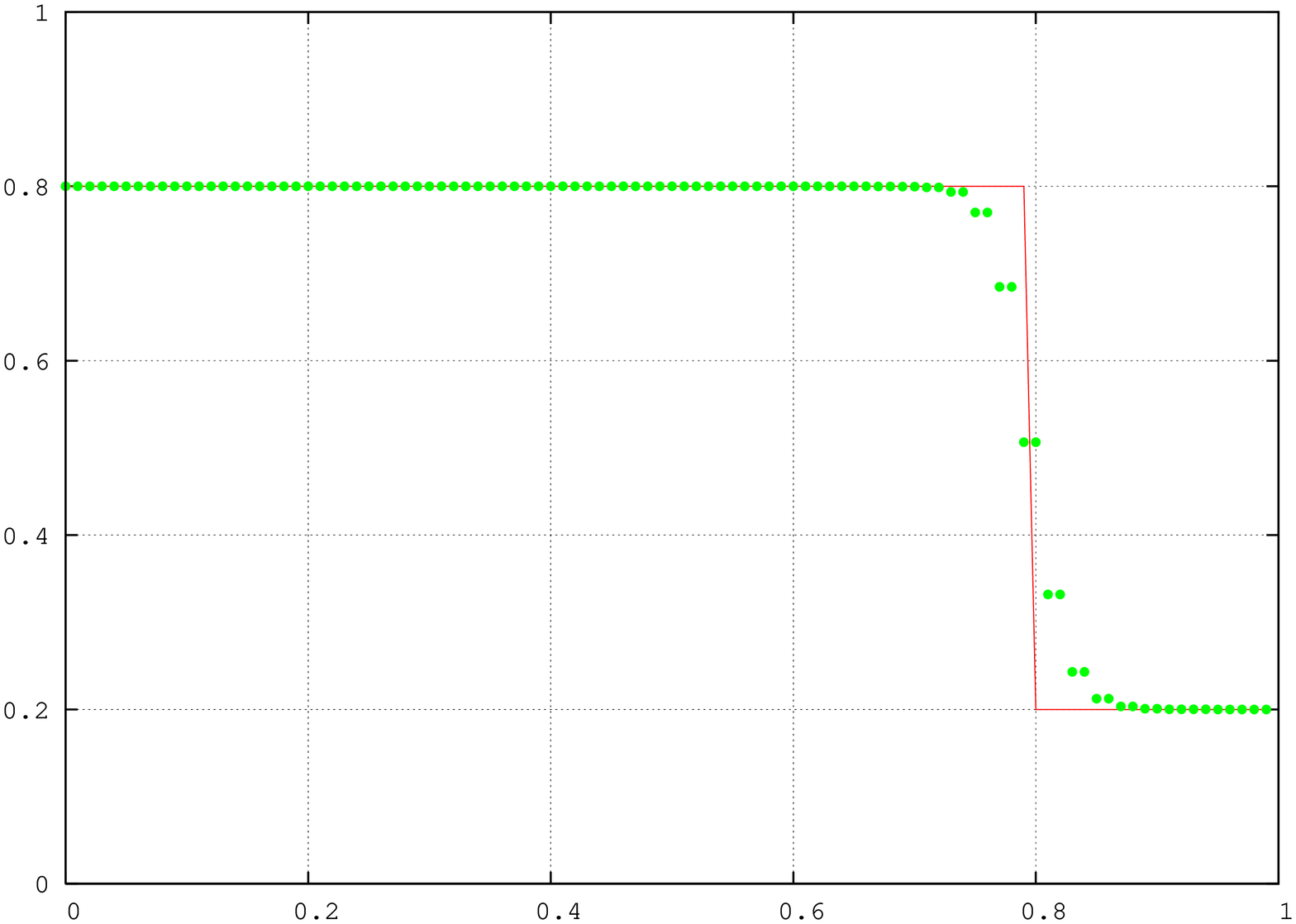}}
\vspace{-1mm} \caption{Lax scheme~(\ref{bya:2}) with Courant
number~$0.9$.}\label{L09}
\end{minipage}
\hfill
\begin{minipage}[b]{7.5cm}
\centerline{\includegraphics[width=5cm,angle=270]{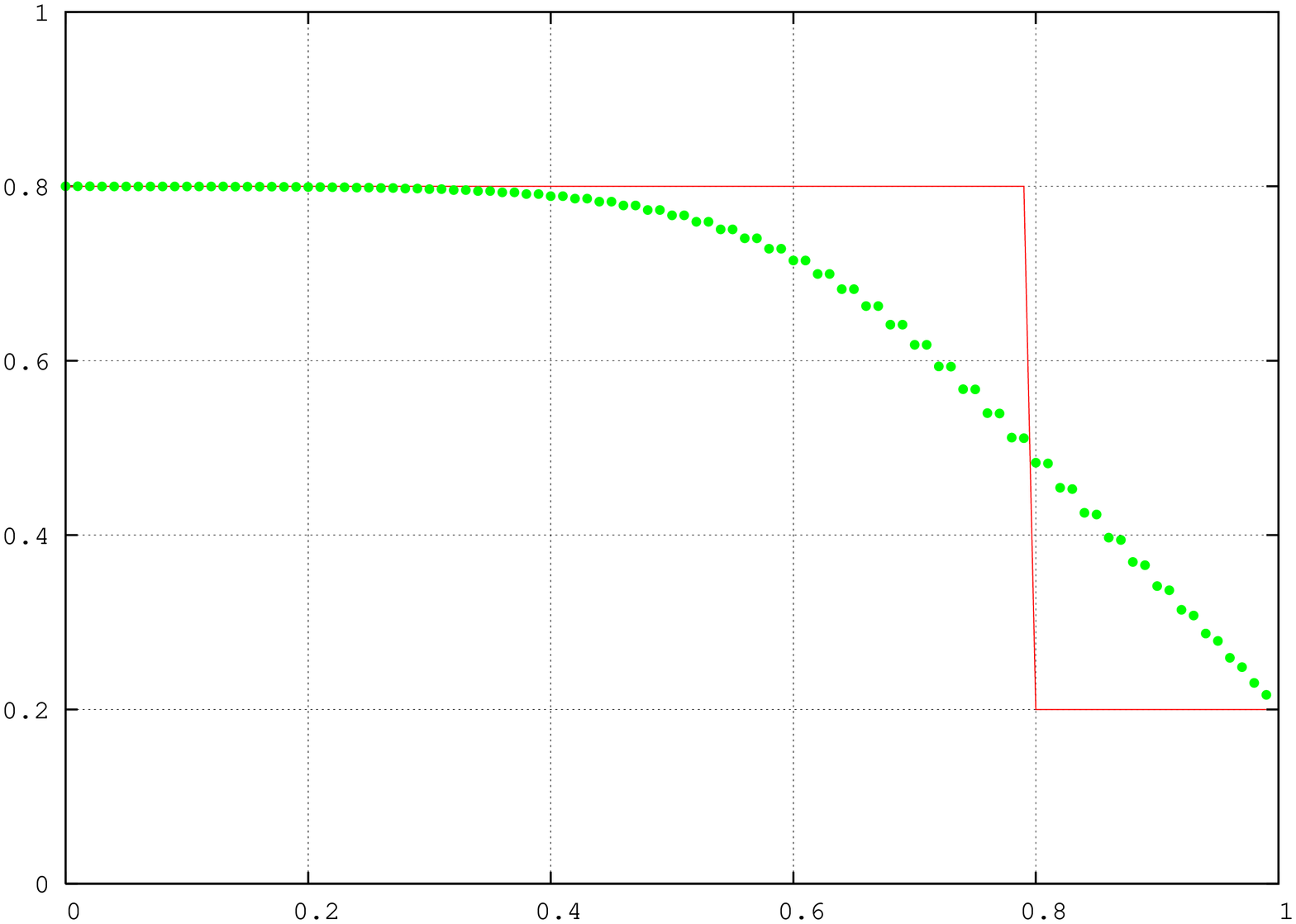}}
\vspace{-1mm} \caption{Lax scheme~(\ref{bya:2}) with Courant
number~$0.1$.}\label{L01}
\end{minipage}
\vspace{-2mm}
\end{figure}

\begin{figure}[t]
\begin{minipage}[b]{7.5cm}
\centerline{\includegraphics[width=5cm,angle=270]{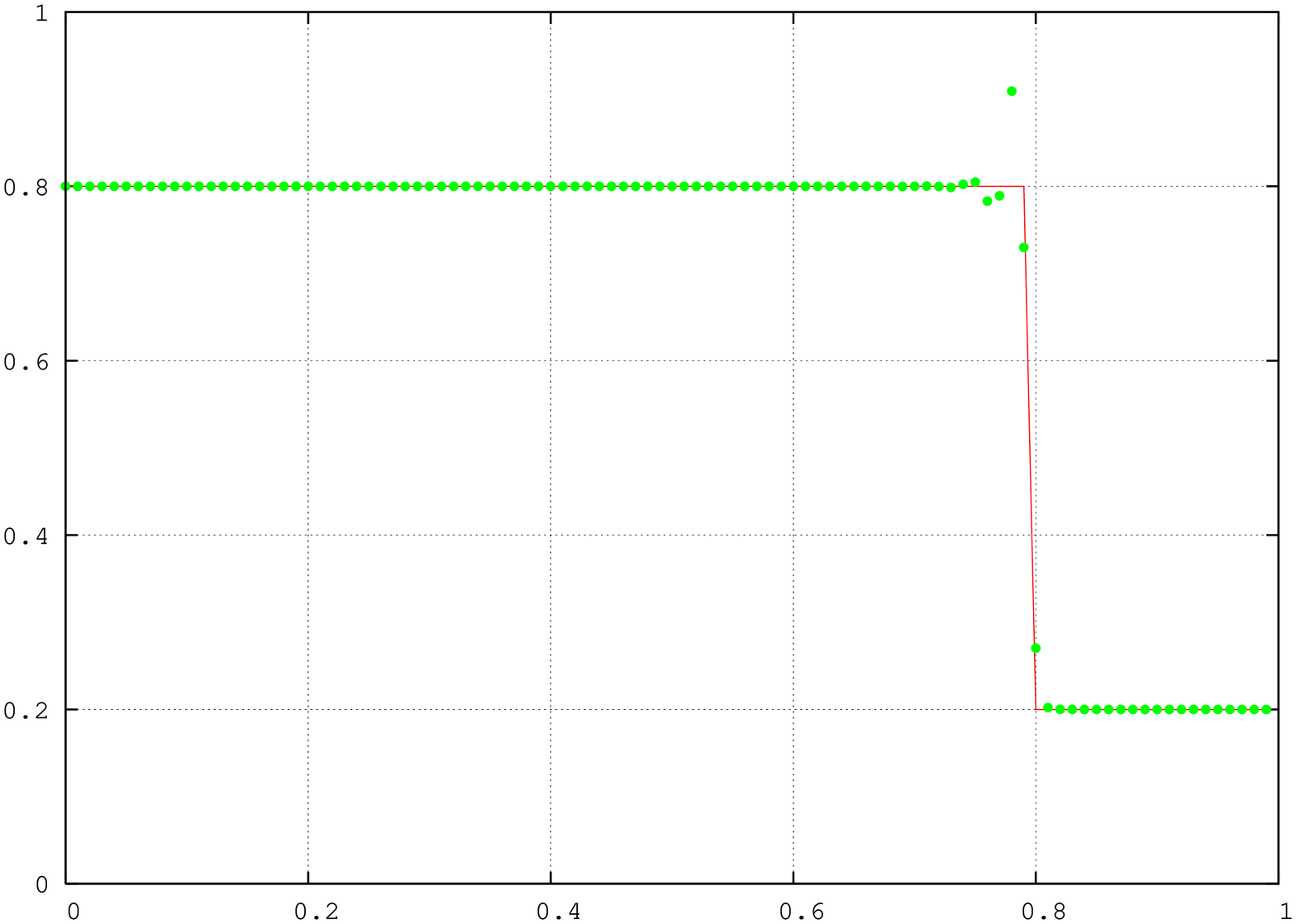}}
\vspace{-1mm} \caption{Lax--Wendrof\/f scheme~(\ref{bya:17}) with
Cou\-rant number~$0.9$.}\label{LW09}
\end{minipage}
\hfill
\begin{minipage}[b]{7.5cm}
\centerline{\includegraphics[width=5cm,angle=270]{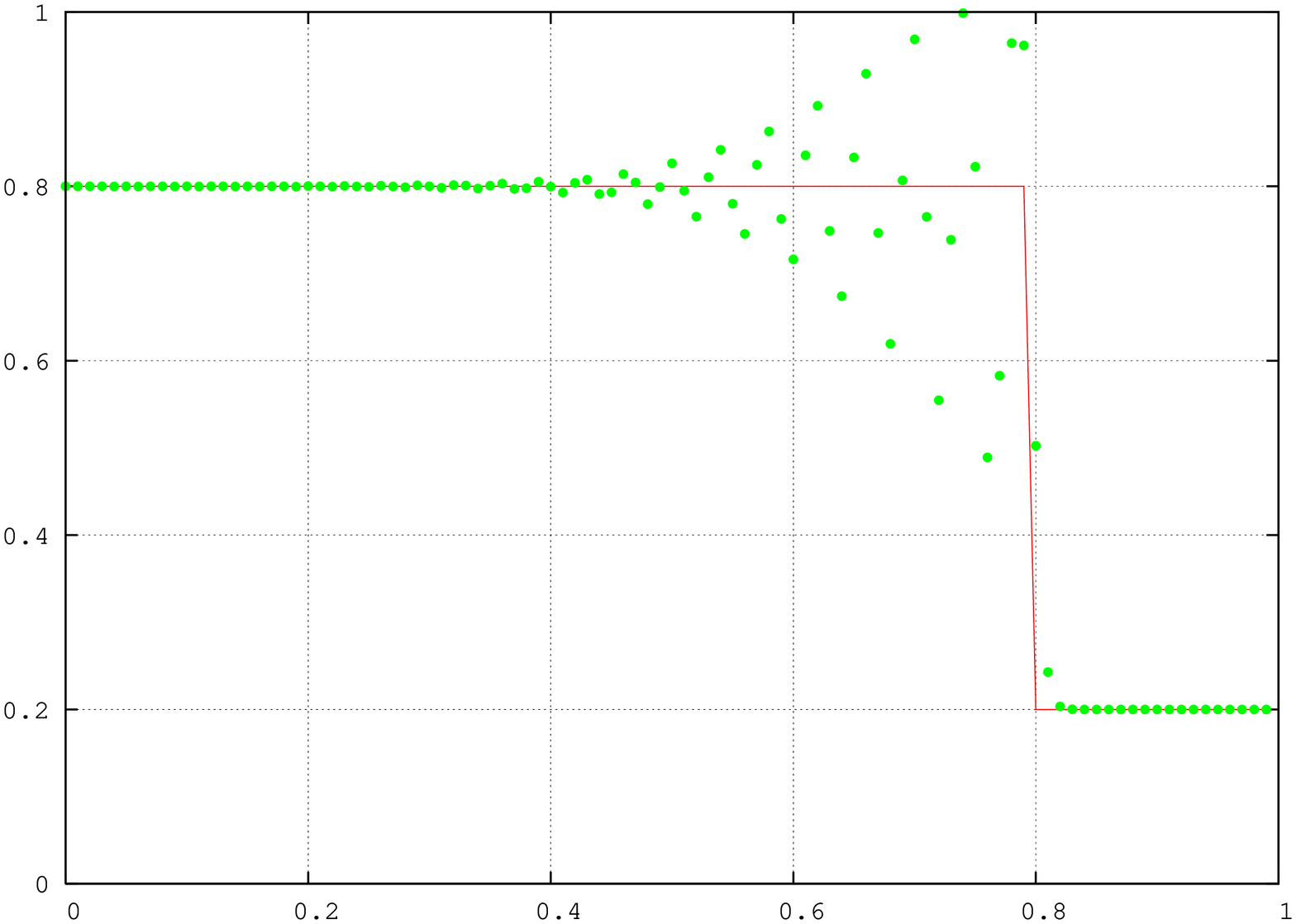}}
\vspace{-1mm} \caption{Lax--Wendrof\/f scheme~(\ref{bya:17}) with
Cou\-rant number $0.1$.}\label{LW01}
\end{minipage}
\vspace{-2mm}
\end{figure}

\begin{figure}[t]
\begin{minipage}[b]{7.5cm}
\centerline{\includegraphics[width=5cm,angle=270]{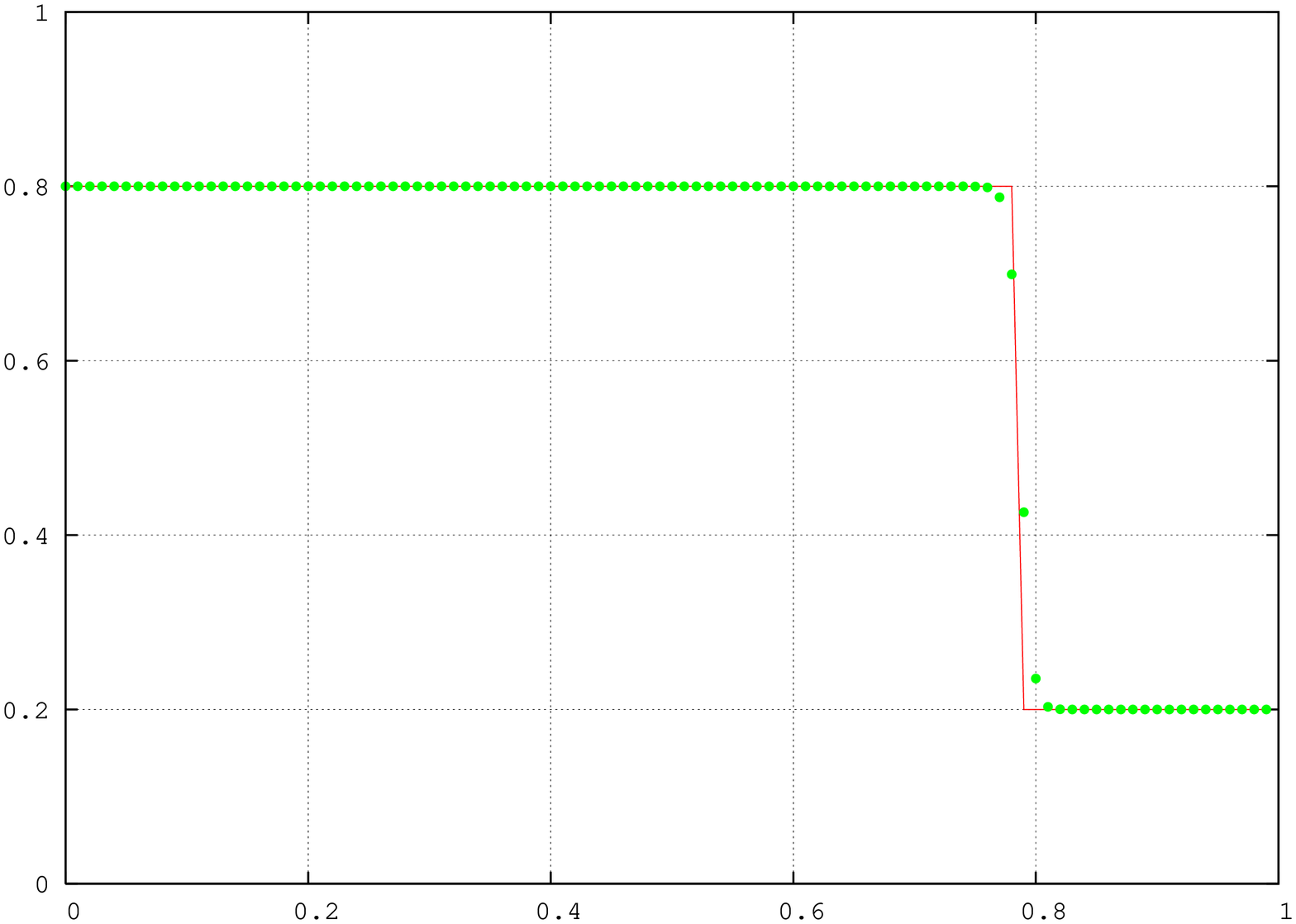}}
\vspace{-1.5mm} \caption{Godunov scheme~(\ref{bya:26}) with
Cou\-rant num\-ber $0.9$.}\label{G09}
\end{minipage}
\hfill
\begin{minipage}[b]{7.5cm}
\centerline{\includegraphics[width=5cm,angle=270]{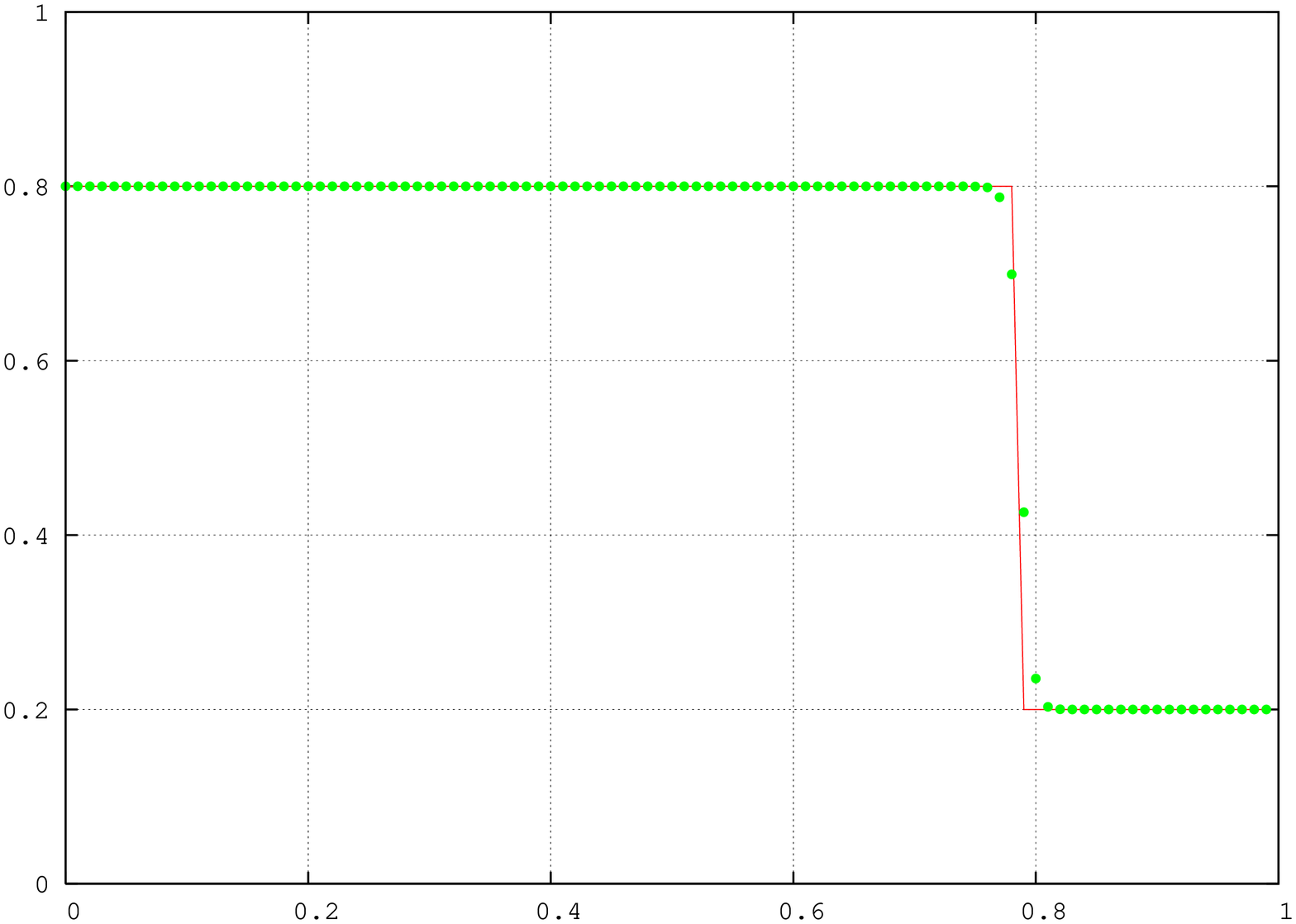}}
\vspace{-1.5mm} \caption{Godunov scheme~(\ref{bya:26}) with
Cou\-rant num\-ber $0.1$.}\label{G01}
\end{minipage}
\vspace{-1.5mm}
\end{figure}

\begin{figure}[t]
\begin{minipage}[t]{7.5cm}
\centerline{\includegraphics[width=6.7cm]{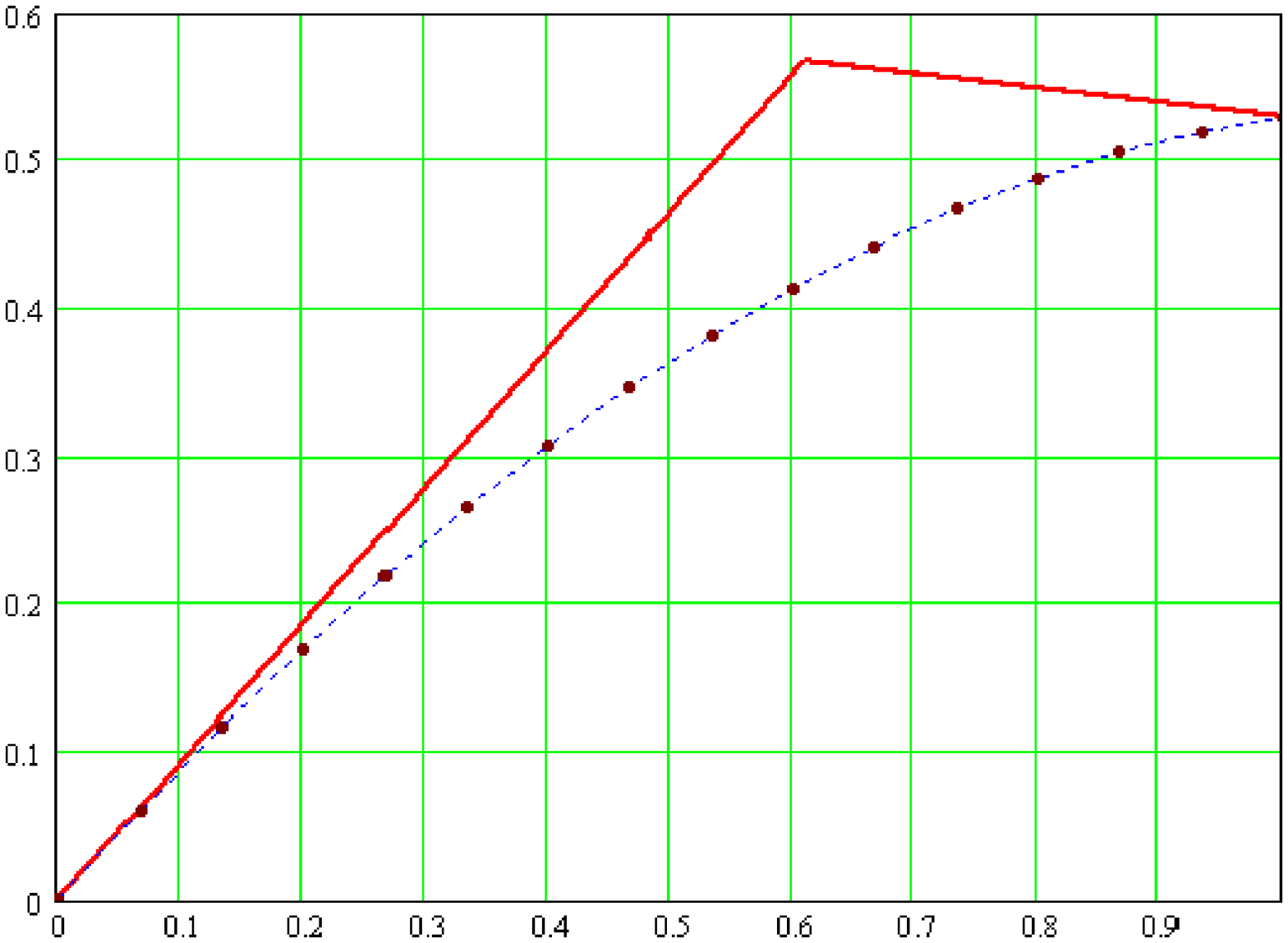}}
\vspace{-1.5mm} \caption{Initial numerical approximation for
equa\-tion~(\ref{FKE_ST}).}\label{IA_FKE}
\end{minipage}
\hfill
\begin{minipage}[t]{7.5cm}
\centerline{\includegraphics[width=6.7cm]{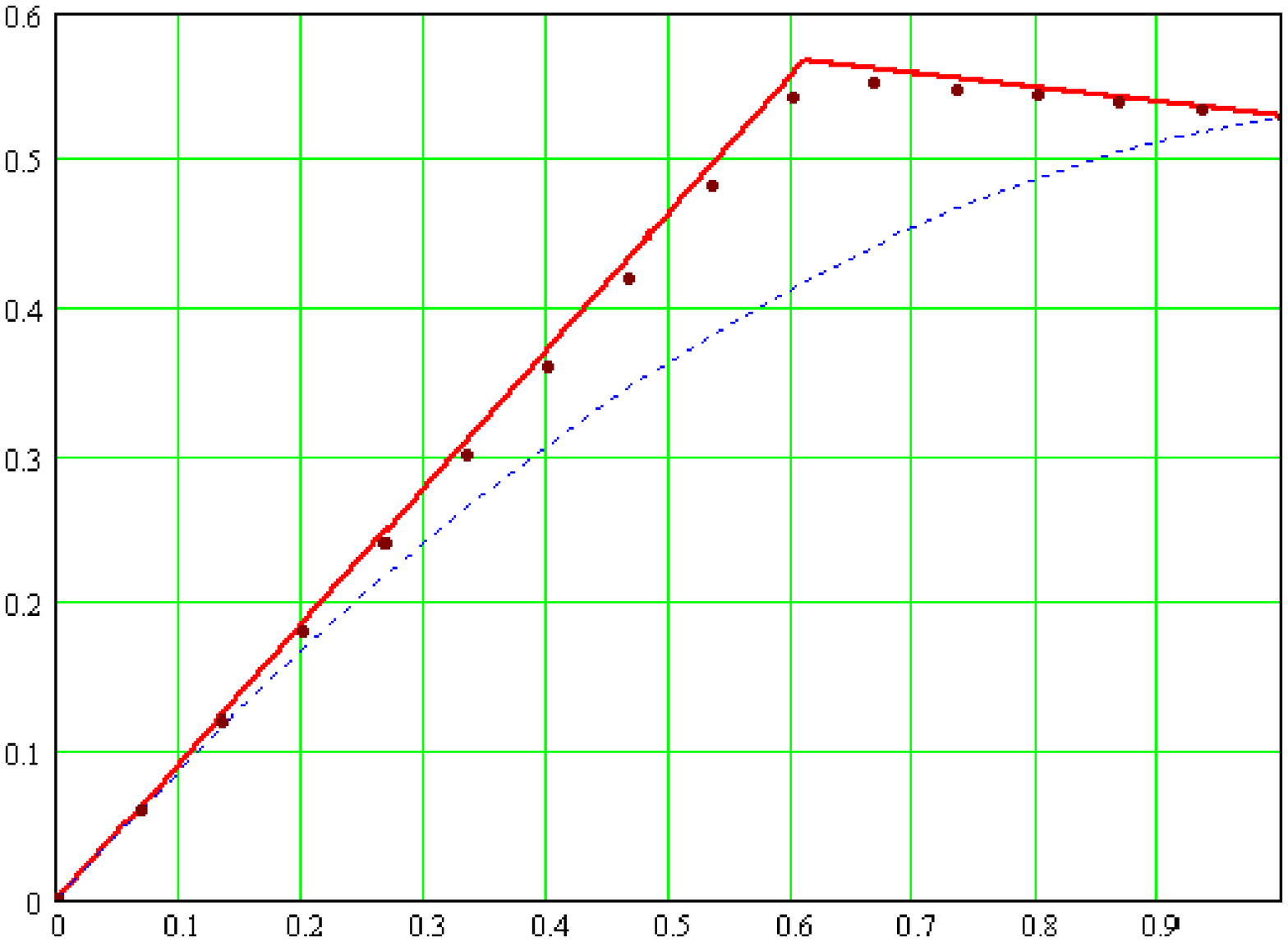}}
\vspace{-1.5mm} \caption{Numerical solution of
equation~(\ref{FKE_ST}).}\label{FA_FKE}
\end{minipage}
\vspace{-4mm}
\end{figure}

As  can be expected (cf.~\cite{GV'96a}), the dispersion ef\/fects
in schemes~(\ref{bya:2})--(\ref{bya:17}) become stronger for the
smaller value of the Courant number (Figs.~\ref{L01}
and~\ref{LW01}). Qualitatively~\cite{Yee'89,Manzini}, the behavior
of Lax scheme~(\ref{bya:2}) in Figs.~\ref{L09}, \ref{L01} is
typical for the classical f\/irst-order schemes when they are
applied to problem~(\ref{BE_nu=0})--(\ref{IC}) whereas the
Lax--Wendrof\/f scheme~(\ref{bya:17}) behaves as the second-oder
method. The Godunov scheme~(\ref{bya:26}), as a shock capturing
one (cf.~\cite{Yee'89,Manzini}), is much better for numerical
description of solution~(\ref{ES}) than the
schemes~(\ref{bya:2})--(\ref{bya:17}), and does not reveal its
sensitivity to the value of the Courant number.

\subsection[Falkowich-Karman equation]{Falkowich--Karman equation}

Now we consider the application of dif\/ference
scheme~(\ref{S_FKE_ST}) to the one-dimensional stationary
transonic f\/low in a channel with a straight density
jump~\cite{J'76}. The exact shock-wave solution of
equation~(\ref{FKE_ST}) at $0\leq x\leq 1$ is shown in
Figs.~\ref{IA_FKE} and~\ref{FA_FKE} by solid red line. Circles
depict the numerical data obtained from dif\/ference
scheme~(\ref{S_FKE_ST}). As an initial approximation, the parabola
was chosen satisfying the following boundary conditions: at the
left, both the function and its derivative are f\/ixed by the
values from the exact solution; at the right, the only function is
bound to the exact solution.

As one can see from~Fig.~\ref{FA_FKE}, scheme~(\ref{S_FKE_ST})
possesses a stable and uniform convergence to the exact shock-wave
solution. Because, by its construction, the scheme is fully
conservative, it does not reveal non-uniqueness of solutions that
is typical for the traditional dif\/ference schemes~\cite{J'76}.

Moreover, the size of the shock transition zone is just one
spatial mesh step that is a consequence of preserving at the
discrete level of all algebraic properties of the initial
PDE~(\ref{FKE}). This is a result of algebraic dif\/ference
elimination provided by the Gr\"obner bases method. Another merit
of scheme~(\ref{S_FKE_ST}) is that it does not involve switches
that is typical for computing transonic f\/low as we already
pointed out in Section~7.

This example shows a principal possibility of constructing
dif\/ference schemes for transonic f\/low without switches and
with the same stencil for both subsonic and supersonic f\/low.

\section{Conclusion}

In the present paper we have shown that the Gr\"obner bases
method, being a universal algorithmic tool for linear dif\/ference
algebra, can be ef\/fectively applied to the construction of
dif\/ferences schemes for linear PDEs with two independent
variables and with rational function coef\/f\/icients. Owing to
the Gr\"obner bases,  this construction is an algorithmic
procedure. It consists in elimi\-na\-tion of partial derivatives
from the system of dif\/ference equations composed from a discrete
version of the original PDEs (on an orthogonal uniform grid) and
numerically approximated integral relations between the unknown
functions and their partial derivatives. As this takes place, the
dif\/ference scheme obtained may depend on the choice of the
integration contour and numerical approximations for integral
relations.

The method is especially ef\/f\/icient when a PDE or a system of
PDEs admits the conservation law form. In this case the
dif\/ference schemes obtained are fully conservative. The
structure of a scheme generated may depend on the choice of
integration contour and numerical integration rules. In so doing,
it is not clear a priori which integration rule leads to a better
scheme.

We also described an ef\/f\/icient algorithm for the construction
of Gr\"obner bases for linear dif\/ference ideals. The algorithm
is based on the concept of Janet-like reductions. Its f\/irst
implementation in Maple is already available, and we used this
implementation in the generation of all linear dif\/ference
schemes presented in the paper.

For classical linear PDEs such as the Laplace equation, the Heat
equation, the Wave equation and the Advection equation our
algorithmic technique leads to the well-known f\/inite
dif\/ference schemes. For Burgers equation we generated several
schemes based on the Lax and Lax--Wendrof\/f methods and computed
their numerical dissipation and dispersion by the dif\/ferential
approximation (modif\/ied equation) method. By example of Burgers
equation we also demonstrate that it is possible to combine the
Godunov method with Gr\"obner bases to derive a shock capturing
scheme.

The non-traditional cubic nonlinear dif\/ference scheme generated
by our dif\/ference elimination method for the Falkowich--Karman
equation describing transonic f\/low in gas dynamics possesses a
number of attractive properties in comparison with traditional
schemes. Among them there are a stable convergence in time to the
exact solution with a one-dimensional shock wave and absence of
switches. It should be noted, however, that due to its cubic
nonlinearity, scheme~(\ref{S_FKE_ST}) has a slower convergence in
comparison with the traditional schemes specially optimized for
numerical simulation of transonic f\/lows in gas dynamics. By this
reason one needs additional research for optimizing nonlinear
schemes obtained by the dif\/ference elimination.

As we already mentioned in the introduction, algorithm {\bf
Gr\"{o}bnerBasis} admits a genera\-li\-zation to
poly\-no\-mial-non\-linear systems of dif\/ference equations
exactly in the same way as the dif\/ferential involutive algorithm
of paper~\cite{Gerdt'99}. In doing so, if every equation in the
initial system is linear with respect to the highest ranking
dif\/ference term and this property of the system is not violated
during its completion to involution, then algorithm {\bf
Gr\"{o}bnerBasis} will work correctly and provide the desirable
output. Such is indeed the case for system~(\ref{DF_of_FKE}). In
the most general case of a dif\/ference system with polynomial
nonlinearity, it can be split into a f\/inite number of subsystems
such that every subsystem can be converted into the Gr\"obner
basis form by applying our algorithm. The underlying splitting
algorithm is a dif\/ference analogue of that described
in~\cite{Thomas}. The latter algorithm is similar to the splitting
algorithm implemented in the library package {\it diffalg} in
Maple.

The above described approach can be also generalized to PDEs with
three and more independent variables. Thus, if PDEs admit the
conservation law form, then one can use multidimensional analogues
of equations (\ref{cons_law}) and (\ref{int_cons_law}) together
with their elementary volume discretization.

\subsection*{Acknowledgements}
We would like to thank the referees for their important remarks
that allowed us to correct the manuscript. We are also grateful to
Daniel Robertz and Viktor Levandovskyy for useful discussions and comments. The
contribution of two authors (V.P.G. and Yu.A.B.) was partially
supported by grants 04-01-00784 and 05-02-17645 from the Russian
Foundation for Basic Research and by grant 2339.2003.2 from the
Ministry of Education and Science of the Russian Federation.

\LastPageEnding

\end{document}